\documentclass{amsart}

\usepackage[leqno]{amsmath}
\usepackage{enumerate, latexsym, mathrsfs, mdwlist}

\swapnumbers
\theoremstyle{definition}
\newtheorem{nul}{}[section]
\newtheorem{dfn}[nul]{Definition}

\newtheorem{cnstr}[nul]{Construction}
\newtheorem{ntn}[nul]{Notation}
\newtheorem{exm}[nul]{Example}

\newtheorem{rec}[nul]{Recollection}
\newtheorem{rmk}[nul]{Remark}

\newtheorem{wrn}[nul]{Warning}
\newtheorem*{dfn*}{Definition}
\newtheorem*{axm*}{Axiom}
\newtheorem*{ntn*}{Notation}
\newtheorem*{exm*}{Example}
\newtheorem*{exr*}{Exercise}
\newtheorem*{int*}{Intuition}
\newtheorem*{qst*}{Question}

\theoremstyle{plain}

\newtheorem{thm}[nul]{Theorem}
\newtheorem{prp}[nul]{Proposition}
\newtheorem{lem}[nul]{Lemma}

\newtheorem{cor}{Corollary}[nul]
\newtheorem*{thm*}{Theorem}
\newtheorem*{prp*}{Proposition}
\newtheorem*{cor*}{Corollary}
\newtheorem*{lem*}{Lemma}
\newtheorem*{cnj*}{Conjecture}

\numberwithin{equation}{nul}

\DeclareMathOperator{\colim}{colim}

\DeclareMathOperator{\Exc}{Exc}

\DeclareMathOperator{\Fun}{Fun}

\DeclareMathOperator{\id}{id}

\DeclareMathOperator{\Map}{Map}
\DeclareMathOperator{\Mor}{Mor}

\DeclareMathOperator{\pr}{pr}

\DeclareMathOperator{\Thy}{Thy}

\newcommand{\BB}{\mathbf{B}}

\newcommand{\D}{\mathrm{D}}

\newcommand{\FF}{\mathbf{F}}

\newcommand{\RR}{\mathbf{R}}
\renewcommand{\SS}{\mathbf{S}}

\newcommand{\VV}{\mathbf{V}}
\newcommand{\WW}{\mathbf{W}}

\newcommand{\ZZ}{\mathbf{Z}}

\newcommand{\Add}{\mathrm{Add}}

\newcommand{\Cat}{\mathbf{Cat}}
\newcommand{\Coh}{\mathbf{Coh}}
\newcommand{\coWald}{\mathbf{coWald}}
\newcommand{\Exact}{\mathbf{Exact}}

\newcommand{\Kan}{\mathbf{Kan}}

\newcommand{\Mod}{\mathbf{Mod}}

\newcommand{\Perf}{\mathbf{Perf}}

\newcommand{\Set}{\mathbf{Set}}
\newcommand{\Sp}{\mathbf{Sp}}

\newcommand{\Wald}{\mathbf{Wald}}

\newcommand{\op}{\mathrm{op}}

\newcommand{\coloneq}{\mathrel{\mathop:}=}

\makeatletter 
\def\revddots{\mathinner{\mkern1mu\raise\p@ 
\vbox{\kern7\p@\hbox{.}}\mkern2mu 
\raise4\p@\hbox{.}\mkern2mu\raise7\p@\hbox{.}\mkern1mu}} 
\makeatother

\usepackage{tikz}
\usetikzlibrary{matrix,arrows,decorations}

\pgfdeclaredecoration{single line}{initial}{\state{initial}[width=\pgfdecoratedpathlength-1sp]{\pgfpathmoveto{\pgfpointorigin}}\state{final}{\pgfpathlineto{\pgfpointorigin}}} 

\newcommand{\fromto}[2]{{#1}\ \tikz[baseline]\draw[>=stealth,->](0,0.5ex)--(0.5,0.5ex);\ {#2}}
\newcommand{\into}[2]{{#1}\ \tikz[baseline]\draw[>=stealth,right hook->](0,0.5ex)--(0.5,0.5ex);\ {#2}}

\newcommand{\cofto}[2]{{#1}\ \tikz[baseline]\draw[>=stealth,>->](0,0.5ex)--(0.5,0.5ex);\ {#2}}
\newcommand{\fibto}[2]{{#1}\ \tikz[baseline]\draw[>=stealth,->>](0,0.5ex)--(0.5,0.5ex);\ {#2}}
\newcommand{\equivto}[2]{{#1}\ \tikz[baseline]\draw[>=stealth,->,font=\scriptsize,inner sep=0.5pt](0,0.5ex)--node[above]{$\sim$}(0.5,0.5ex);\ {#2}}
\newcommand{\goesto}[2]{{#1}\ \tikz[baseline]\draw[|->](0,0.5ex)--(0.5,0.5ex);\ {#2}}

\renewcommand{\to}{\ \tikz[baseline]\draw[>=stealth,->](0,0.5ex)--(0.5,0.5ex);\ }

\title{On exact $\infty$-categories and the {T}heorem of the {H}eart}

\author{Clark Barwick}
\email{clarkbar@gmail.com}
\address{Massachusetts Institute of Technology, Department of Mathematics, Building 2, 77 Massachusetts Avenue, Cambridge, MA 02139-4307, USA}

\begin{document}

\begin{abstract} The new homotopy theory of \emph{exact $\infty$-categories} is introduced and employed to prove a Theorem of the Heart for algebraic $K$-theory (in the sense of Waldhausen). This implies a new compatibility between Waldhausen $K$-theory and Neeman $K$-theory. Additionally, it provides a new proof of the D\'evissage and Localization Theorems of Blumberg--Mandell, new models for the $G$-theory of schemes, and a proof of the invariance of $G$-theory under derived nil-thickenings.
\end{abstract}

\maketitle

\setcounter{tocdepth}{2}
\tableofcontents

\setcounter{section}{-1}

\section{Introduction} In this paper, we prove (Th. \ref{thm:thmofheart}) that a stable homotopy theory whose triangulated homotopy category admits a bounded t-structure has the same algebraic $K$-theory (in the sense of Waldhausen) as that of its heart. This is the \emph{Theorem of the Heart} of the title. This result, which apparently has been expected by some experts for some time, has nevertheless gone unproved.

This theorem does, however, have an important predecessor. For 20 years Amnon Neeman has advanced the \emph{algebraic $K$-theory of triangulated categories}  \cite{MR1491990,MR1604910,MR1656552,MR1724625,MR1793672,MR1798824,MR1798828,MR1828612,MR2181838} as a way of extracting $K$-theoretic data directly from the triangulated homotopy category of a stable homotopy theory. As an approximation to Waldhausen $K$-theory, this form of $K$-theory has well-documented limitations: a beautiful example of Marco Schlichting \cite{MR1930883} shows that Waldhausen $K$-theory can distinguish stable homotopy theories with equivalent triangulated homotopy categories. Nevertheless, the most impressive advance in the algebraic $K$-theory of triangulated categories is Neeman's \emph{Theorem of the Heart} \cite{MR1724625,MR1793672,MR1828612}, which expresses an equivalence between the Neeman $K$-theory of a triangulated category $\mathscr{T}$ equipped with a bounded t-structure and the Quillen $K$-theory of its heart $\mathscr{T}^{\heartsuit}$.

Neeman's proof of his Theorem of the Heart is lengthy and difficult to read, so much so that it even generated a small controversy (see Neeman's discussion in \cite[pp. 347--353]{MR1491990}). The proof of our Theorem of the Heart, by contrast, is mercifully short, conceptually appealing, and logically independent of Neeman's. Consequently we regard Th. \ref{thm:thmofheart} and its proof as a conclusive answer to Problem 76 of his survey \cite{MR2181838}.

To prove our result, we introduce a natural homotopy-theoretic generalization of Quillen's notion of an exact category, which we call an \emph{exact $\infty$-category} (Df. \ref{dfn:biWald}). Because this notion involves a compatibility between certain homotopy limits and certain homotopy colimits, it is virtually impossible (or at the very least hideously inconvenient) to express in the more classical language of categories-with-weak-equivalences. Therefore we have to employ concepts from higher category theory --- in particular, the theory of \emph{Waldhausen $\infty$-categories}, whose theory we studied in pitiless detail in \cite{K1}.

The key idea from \cite{K1} is that algebraic $K$-theory is a homology theory for $\infty$-categories. In fact, algebraic $K$-theory is the analogue of stable homotopy theory in this context. The behavior of these categorified homology theories under duality is the key phenomenon that makes our proof of the Theorem of the Heart work. More precisely, when algebraic $K$-theory is restricted to exact $\infty$-categories, it enjoys a \emph{self-duality} (Cor. \ref{cor:Kdual}). This self-duality is then used in conjunction with our $\infty$-categorical Fibration Theorem \cite[Pr. 10.12]{K1} to prove the following.
\begin{thm*}[Heart] If $\mathscr{A}$ is a stable $\infty$-category equipped with a bounded t-structure, then the inclusion of the heart $\mathscr{A}^{\heartsuit}\subset\mathscr{A}$ induces a $K$-theory weak equivalence
\begin{equation*}
K(\mathscr{A}^{\heartsuit})\simeq K(\mathscr{A}).
\end{equation*}
\end{thm*}
\noindent This result is one of the very few general statements in algebraic $K$-theory that is capable of providing $K$-theory equivalences that do not arise from equivalences of the $\infty$-categories themselves. (The only other example of such a general result we know of this kind is Quillen's D\'evissage Theorem.)

The full strength of the conceptual apparatus constructed here and in \cite{K1} is necessary for this proof to work. In view of Schlichting's ``no-go theorem'' \cite[Pr. 2.2]{MR1930883}, our use of the Fibration Theorem makes it impossible for a proof at all similar to the one presented here to be adapted to the context of triangulated categories. On the other hand, there are other proposed versions of algebraic $K$-theory for $\infty$-categories (most notably that of Blumberg, Gepner, and Tabuada \cite{BGT}) that restrict attention to stable $\infty$-categories or the like. These versions of $K$-theory just won't help for this problem: in fact, it isn't possible even to \emph{express} the relevant cases of self-duality with this or any other form of $K$-theory that splits arbitrary cofiber sequences.

Let us underscore that this is not a new proof of an old theorem. Schlichting's example shows that there is no \emph{a priori} reason to expect Neeman's Theorem of the Heart to say anything about Waldhausen $K$-theory. Nevertheless, our main result does yield a comparison between Neeman's $K$-theory and Waldhausen $K$-theory. Indeed, the conjunction of Neeman's Theorem of the Heart and Th. \ref{thm:thmofheart} implies that the Waldhausen $K$-theory of a stable $\infty$-category $\mathscr{A}$ agrees with the Neeman $K$-theory of its triangulated homotopy category $\mathscr{T}=h\mathscr{A}$ (the variant denoted $K({}^w\!\mathscr{T})$ in \cite{MR2181838}), whenever the latter admits a bounded t-structure (Cor. \ref{cor:Neemanconj}). This verifies a conjecture of Neeman \cite[Conj. A.5]{MR1793672} for such stable homotopy theories.

This paper ends with a discussion of some immediate corollaries of the main theorem, which we include mostly as proof of concept.
\begin{enumerate}
\item[---] We give new models for the $G$-theory of schemes in terms of the perverse coherent sheaves of Arinkin, Bezrukavnikov, and Deligne (Ex. \ref{exm:perv}).
\item[---] We also give a new, short proof of the D\'evissage and Localization Theorems of Blumberg--Mandell \cite{BM} (Pr. \ref{prp:dev} and Th. \ref{thm:BM}), which immediately yields a host of useful fiber sequences in the algebraic $K$-theory of ring spectra (Ex. \ref{exm:BMfiberseqs}). More interesting examples can be found in our paper with Tyler Lawson \cite{BLreg}.
\item[---] Finally, we show that the $G$-theory of spectral Deligne--Mumford stacks (in the sense of Lurie) is invariant under derived thickenings (Pr. \ref{prp:Gthynothingnew}).
\end{enumerate}

\subsection*{Acknowledgments} I thank Andrew Blumberg and Mike Mandell for very helpful conversations about this paper. I also thank Dustin Clausen for helpful remarks about a previous version of this paper. I am very grateful to Benjamin Antieau for noticing that overzealous revision botched the proof of the main theorem in a earlier version of this paper. Finally, I thank the anonymous referee for a genuinely helpful report.

\section{Preliminaries} We use higher categories systematically in this paper. In particular, we are interested in $\infty$-categories whose $i$-morphisms for $i\geq 2$ are all invertible. These days, it is fashionable to call these \emph{$(\infty,1)$-categories} or even just \emph{$\infty$-categories}. There are very many models for the homotopy theory of $\infty$-categories in this sense, and they are all equivalent in an essentially unique fashion, up to orientation --- see To\"en \cite{Toen} or Lurie \cite{G} or Barwick--Schommer-Pries \cite{BSP}.

In this paper, we employ the homotopy theory of \emph{quasicategories} developed by Joyal \cite{Joyal, Joyal08} and then further by Lurie \cite{HTT}. These are simplicial sets $C$ in which any inner horn $x\colon\fromto{\Lambda^m_k}{C}$ ($m\geq2$ and $1\leq k\leq m-1$) admits a filler $\overline{x}\colon\fromto{\Delta^m}{C}$. When we use the phrase ``$\infty$-categories'' in this text, we will be referring to these.

One point that is perhaps not so obvious is the notion of a \emph{subcategory} of an $\infty$-category.
\begin{rec}[\protect{\cite[\S 1.2.11]{HTT}}]\label{rec:subcats} A \emph{subcategory} of an $\infty$-category $A$ is a simplicial subset $A'\subset A$ such that for some subcategory $(hA)'$ of the homotopy category $hA$, the square
\begin{equation*}
\begin{tikzpicture} 
\matrix(m)[matrix of math nodes, 
row sep=4ex, column sep=4ex, 
text height=1.5ex, text depth=0.25ex] 
{A'&A\\ 
N(hA)'&N(hA)\\}; 
\path[>=stealth,->,font=\scriptsize] 
(m-1-1) edge[right hook->] (m-1-2) 
edge (m-2-1) 
(m-1-2) edge (m-2-2) 
(m-2-1) edge[right hook->] (m-2-2); 
\end{tikzpicture}
\end{equation*}
is a pullback diagram of simplicial sets. In particular, note that a subcategory of an $\infty$-category is uniquely specified by specifying a subcategory of its homotopy category. Note also that any inclusion $\into{A'}{A}$ of a subcategory is an inner fibration \cite[Df. 2.0.0.3, Pr. 2.3.1.5]{HTT}.

We will say that $A'\subset A$ is a \emph{full subcategory} if $(hA)'\subset hA$ is a full subcategory. In this case, $A'$ is uniquely determined by the set $A'_0$ of vertices of $A'$, and we say that $A'$ is \emph{spanned by} the set $A'_0$.

We will say that $A'$ is \emph{stable under equivalences} if the subcategory $(hA)'\subset hA$ above can be chosen to be stable under isomorphisms. Note that any inclusion $\into{A'}{A}$ of a subcategory that is stable under equivalences is a categorical fibration, i.e., a fibration for the Joyal model structure \cite[Cor. 2.4.6.5]{HTT}.
\end{rec}

The natural inputs for algebraic $K$-theory are what we call Waldhausen $\infty$-categories:
\begin{rec}[\protect{\cite[Df. 2.7]{K1}}]\label{dfn:Waldinftycat} A \emph{Waldhausen $\infty$-category} $(\mathscr{C},\mathscr{C}_{\dag})$ consists of an $\infty$-category $\mathscr{C}$ equipped with a subcategory $\mathscr{C}_{\dag}\subset\mathscr{C}$ that contains all the equivalences. A morphism of $\mathscr{C}_{\dag}$ will be said to be \emph{ingressive} or \emph{a cofibration}. These data are then required to satisfy the following conditions.
\begin{enumerate}[(\ref{dfn:Waldinftycat}.1)]
\item The $\infty$-category $\mathscr{C}$ contains a \emph{zero object} --- i.e., an object that is both initial and terminal \cite[Df. 1.2.12.1 and Rk. 1.2.12.6]{HTT}.
\item For any zero object $0$, any morphism $\fromto{0}{X}$ is ingressive.
\item Pushouts \cite[\S 4.4.2]{HTT} of cofibrations exist and are cofibrations.
\end{enumerate}

A functor $\fromto{\mathscr{C}}{\mathscr{D}}$ between Waldhausen $\infty$-categories is said to be \emph{exact} if it carries cofibrations to cofibrations and preserves both zero objects and pushouts of cofibrations.
\end{rec}

The examples one may have in mind here includes the nerve of an ordinary exact category in the sense of Quillen (in which the ingressive morphisms are the admissible monomorphisms), the nerve of a category with cofibrations in the sense of Waldhausen \cite{} (in which the ingressive morphisms are the cofibrations), any $\infty$-category with a zero object and all finite colimits (in which any morphism is ingressive).

In a sense, the defining property of algebraic $K$-theory is that it splits cofiber sequences. We'll discuss this point in more detail later. For now, let us explain what cofiber sequences \emph{are}.
\begin{dfn} In a Waldhausen $\infty$-category, a \emph{cofiber sequence} is a pushout square
\begin{equation*}
\begin{tikzpicture} 
\matrix(m)[matrix of math nodes, 
row sep=4ex, column sep=4ex, 
text height=1.5ex, text depth=0.25ex] 
{X'&X\\ 
0&X''\\}; 
\path[>=stealth,->,font=\scriptsize] 
(m-1-1) edge[>->] (m-1-2) 
edge (m-2-1) 
(m-1-2) edge (m-2-2) 
(m-2-1) edge[>->] (m-2-2); 
\end{tikzpicture}
\end{equation*}
in which $\cofto{X'}{X}$ is ingressive and $0$ is a zero object. We call $\fromto{X}{X''}$ the \emph{cofiber} of the cofibration $\cofto{X'}{X}$.
\end{dfn}

We also have the dual notion of a coWaldhausen $\infty$-category.
\begin{rec}[\protect{\cite[Df. 2.16]{K1}}]\label{dfn:coWaldinftycat} A \emph{coWaldhausen $\infty$-category} $(\mathscr{C},\mathscr{C}^{\dag})$ consists of an $\infty$-category $\mathscr{C}$ equipped with a subcategory $\mathscr{C}^{\dag}\subset\mathscr{C}$ that contains all the equivalences. A morphism of $\mathscr{C}^{\dag}$ will be said to be \emph{egressive} or \emph{a fibration}. These data are then required to satisfy the following conditions.
\begin{enumerate}[(\ref{dfn:coWaldinftycat}.1)]
\item The $\infty$-category $\mathscr{C}$ contains a zero object.
\item For any zero object $0$, any morphism $\fromto{X}{0}$ is egressive.
\item Pullbacks \cite[\S 4.4.2]{HTT} of fibrations exist and are fibrations.
\end{enumerate}

A functor $\fromto{\mathscr{C}}{\mathscr{D}}$ between coWaldhausen $\infty$-categories is said to be \emph{exact} if it carries fibrations to fibrations and preserves both zero objects and pullbacks of fibrations.
\end{rec}

In other words, coWaldhausen $\infty$-categories are precisely the opposites of Waldhausen $\infty$-categories. In fact, Waldhausen and coWaldhausen $\infty$-categories and exact functors organize themselves into $\infty$-categories $\Wald_{\infty}$ and $\coWald_{\infty}$, and the formation of the opposite $\infty$-category restricts to an equivalence of $\infty$-categories
\begin{equation*}
\equivto{\Wald_{\infty}}{\coWald_{\infty}}
\end{equation*}
\cite[Nt. 2.13, Nt. 2.17, and Pr. 2.18]{K1}.

\begin{dfn} In a coWaldhausen $\infty$-category, a \emph{fiber sequence} is a pullback square
\begin{equation*}
\begin{tikzpicture} 
\matrix(m)[matrix of math nodes, 
row sep=4ex, column sep=4ex, 
text height=1.5ex, text depth=0.25ex] 
{X'&X\\ 
0&X''\\}; 
\path[>=stealth,->,font=\scriptsize] 
(m-1-1) edge (m-1-2) 
edge[->>] (m-2-1) 
(m-1-2) edge[->>] (m-2-2) 
(m-2-1) edge (m-2-2); 
\end{tikzpicture}
\end{equation*}
in which $\fibto{X}{X''}$ is egressive and $0$ is a zero object. We call $\fromto{X'}{X}$ the \emph{fiber} of the fibration $\fibto{X}{X''}$.
\end{dfn}


\section{Additive $\infty$-categories} In effect, an exact $\infty$-category will be a Waldhausen $\infty$-category $\mathscr{C}$ that is also a coWaldhausen $\infty$-category, with two  additional properties: first, fiber sequences and cofiber sequences in $\mathscr{C}$ must coincide; and second, $\mathscr{C}$ must be \emph{additive}.

This notion of additivity is similar to the notion of additivity for ordinary categories. Recall from \cite[Df. 4.5]{K1} the following.
\begin{rec} \label{item:directsums} Suppose $C$ is an $\infty$-category. Then $C$ is said to \emph{admit finite direct sums} if the following conditions hold.
\begin{enumerate}[(\ref{item:directsums}.1)]
\item The $\infty$-category $C$ admits a zero object.
\item The $\infty$-category $C$ has all finite products and coproducts.
\item For any finite set $I$ and any $I$-tuple $(X_i)_{i\in I}$ of objects of $C$, the map
\begin{equation*}
\fromto{\coprod X_I}{\prod X_I}
\end{equation*}
in $hC$ --- given by the maps $\phi_{ij}\colon\fromto{X_i}{X_j}$, where $\phi_{ij}$ is zero unless $i=j$, in which case it is the identity --- is an isomorphism.
\end{enumerate}
If $C$ admits finite direct sums, then for any finite set $I$ and any $I$-tuple $(X_i)_{i\in I}$ of objects of $C$, we denote by $\bigoplus X_I$ the product (or, equivalently, the coproduct) of the $X_i$.
\end{rec}

Suppose $C$ an $\infty$-category that admits direct sums. Then the homotopy category $hC$ is easily seen to admit direct sums as well. Moreover, the mapping spaces in $C$ admit the natural structure of a homotopy-commutative $H$-space: for any morphisms $f,g\in\Map_{C}(X,Y)$, one may define $f+g\in\Map_{C}(X,Y)$ as the composite
\begin{equation*}
X\ \tikz[baseline]\draw[>=stealth,->,font=\scriptsize](0,0.5ex)--node[above]{$\Delta$}(0.5,0.5ex);\ X\oplus X\ \tikz[baseline]\draw[>=stealth,->,font=\scriptsize](0,0.5ex)--node[above]{$f\oplus g$}(0.85,0.5ex);\ Y\oplus Y\ \tikz[baseline]\draw[>=stealth,->,font=\scriptsize](0,0.5ex)--node[above]{$\nabla$}(0.5,0.5ex);\ Y.
\end{equation*}

\begin{dfn} An $\infty$-category $C$ that admits finite direct sums will be said to be \emph{additive} if its homotopy category $hC$ is additive; i.e., if for any two objects $X$ and $Y$, the monoid $\pi_0\Map_C(X,Y)$ is a group.
\end{dfn}

\begin{rmk} An $\infty$-category $C$ with direct sums is additive just in case, for any objects $X$ and $Y$, the \emph{shear map}
\begin{equation*}
\fromto{\Map_C(X,Y)\times\Map_C(X,Y)}{\Map_C(X,Y)\times\Map_C(X,Y)}
\end{equation*}
in the homotopy category of Kan simplicial sets that is given informally by the assignment $\goesto{(f,g)}{(f,f+g)}$ is an isomorphism. Note in particular that additivity is a condition, not additional structure.
\end{rmk} 

\begin{exm} Clearly the nerve of any ordinary additive category is an additive $\infty$-category. Similarly, any stable $\infty$-category is additive.
\end{exm}


\section{Exact $\infty$-categories} Now we are ready to define exact $\infty$-categories.

\begin{dfn}\label{dfn:biWald} Suppose $\mathscr{C}$ an $\infty$-category, and suppose $\mathscr{C}_{\dag}$ and $\mathscr{C}^{\dag}$ two subcategories of $\mathscr{C}$. Call a morphism of $\mathscr{C}_{\dag}$ \emph{ingressive} or \emph{a cofibration}, and call a morphism of $\mathscr{C}^{\dag}$ \emph{egressive} or \emph{a fibration}.
\begin{enumerate}[(\ref{dfn:biWald}.1)]
\item A pullback square
\begin{equation*}
\begin{tikzpicture}[baseline]
\matrix(m)[matrix of math nodes, 
row sep=4ex, column sep=4ex, 
text height=1.5ex, text depth=0.25ex] 
{X&Y\\
X'&Y',\\}; 
\path[>=stealth,->,font=\scriptsize] 
(m-1-1) edge (m-1-2) 
edge (m-2-1) 
(m-1-2) edge[->>] (m-2-2) 
(m-2-1) edge[>->] (m-2-2); 
\end{tikzpicture}
\end{equation*}
is said to be \emph{ambigressive} if $\cofto{X'}{Y'}$ is ingressive and $\fibto{Y}{Y'}$ is egressive. Dually, a pushout square
\begin{equation*}
\begin{tikzpicture}[baseline]
\matrix(m)[matrix of math nodes, 
row sep=4ex, column sep=4ex, 
text height=1.5ex, text depth=0.25ex] 
{X&Y\\
X'&Y',\\}; 
\path[>=stealth,->,font=\scriptsize] 
(m-1-1) edge[>->] (m-1-2) 
edge[->>] (m-2-1) 
(m-1-2) edge (m-2-2) 
(m-2-1) edge (m-2-2); 
\end{tikzpicture}
\end{equation*}
is said to be \emph{ambigressive} if $\cofto{X}{Y}$ is ingressive and $\fibto{X}{X'}$ is egressive.
\item\label{item:exactinfty} The triple $(\mathscr{C},\mathscr{C}_{\dag},\mathscr{C}^{\dag})$ is said to be an \emph{exact $\infty$-category} if it satisfies the following conditions.
\begin{enumerate}[(\ref{dfn:biWald}.\ref{item:exactinfty}.1)]
\item The underlying $\infty$-category $\mathscr{C}$ is additive.
\item The pair $(\mathscr{C},\mathscr{C}_{\dag})$ is a Waldhausen $\infty$-category.
\item The pair $(\mathscr{C},\mathscr{C}^{\dag})$ is a coWaldhausen $\infty$-category.
\item A square in $\mathscr{C}$ is an ambigressive pullback if and only if it is an ambigressive pushout.
\end{enumerate}
\item In an exact $\infty$-category, an \emph{exact sequence} is a cofiber/fiber sequence
\begin{equation*}
\begin{tikzpicture} 
\matrix(m)[matrix of math nodes, 
row sep=4ex, column sep=4ex, 
text height=1.5ex, text depth=0.25ex] 
{X'&X\\ 
0&X'';\\}; 
\path[>=stealth,->,font=\scriptsize] 
(m-1-1) edge[>->] (m-1-2) 
edge[->>] (m-2-1) 
(m-1-2) edge[->>] (m-2-2) 
(m-2-1) edge[>->] (m-2-2); 
\end{tikzpicture}
\end{equation*}
we will abuse notation by writing
\begin{equation*}
X'\ \tikz[baseline]\draw[>=stealth,>->](0,0.5ex)--(0.5,0.5ex);\ X\ \tikz[baseline]\draw[>=stealth,->>](0,0.5ex)--(0.5,0.5ex);\ X''
\end{equation*}
for this square.
\end{enumerate}
\end{dfn}

\begin{rmk} Note that in an exact $\infty$-category, a morphism of an exact $\infty$-category is egressive just in case it appears as the cofiber of an ingressive morphism, and, dually, a morphism of an exact $\infty$-category is ingressive just in case it appears as the fiber of an egressive morphism. Indeed, any cofiber of an ingressive morphism is egressive, and any egressive morphism is equivalent to the cofiber of its fiber. This proves the first statement; the second is dual. Consequently, the class of cofibrations in an exact $\infty$-category specifies the class of fibrations, and vice versa.
\end{rmk}

\begin{exm}\label{exm:exactinftycats}\begin{enumerate}[(\ref{exm:exactinftycats}.1)]
\item The nerve $NC$ of an ordinary category $C$ can be endowed with a triple structure yielding an exact $\infty$-category if and only if $C$ is an ordinary exact category, in the sense of Quillen, wherein the admissible monomorphisms are exactly the cofibrations, and the admissible epimorphisms are exactly the fibrations. To prove this, one may observe that the ``minimal'' axioms of Keller \cite[App. A]{MR1052551} simply \emph{are} the axioms listed above.
\item At the other extreme, any stable $\infty$-category is an exact $\infty$-category in which all morphisms are both egressive and ingressive, and, conversely, any $\infty$-category that can be regarded as an exact category with the \emph{maximal} triple structure (in which any morphism is both ingressive and egressive) is a stable $\infty$-category.
\end{enumerate}
\end{exm}

\begin{exm}\label{exm:tstructexacts} We may interpolate between these two extremes. Suppose $\mathscr{A}$ a stable $\infty$-category equipped with a t-structure, and suppose $a,b\in\ZZ$.
\begin{enumerate}[(\ref{exm:tstructexacts}.1)]
\item The $\infty$-category $\mathscr{A}_{[a,+\infty)}\coloneq\mathscr{A}_{\geq a}$ admits an exact $\infty$-category structure, in which every morphism is ingressive, but a morphism $\fromto{Y}{Z}$ is egressive just in case the induced morphism $\fromto{\pi_aX}{\pi_aY}$ is an epimorphism of $\mathscr{A}^{\heartsuit}$.
\item Dually, the $\infty$-category $\mathscr{A}_{(-\infty,b]}\coloneq\mathscr{A}_{\leq b}$ admits an exact $\infty$-category structure, in which every morphism is egressive, but a morphism $\fromto{X}{Y}$ is ingressive just in case the induced morphism $\fromto{\pi_bX}{\pi_bY}$ is a monomorphism of $\mathscr{A}^{\heartsuit}$.
\item We may intersect these subcategories to obtain the full subcategory
\begin{equation*}
\mathscr{A}_{[a,b]}\coloneq\mathscr{A}_{\geq a}\cap\mathscr{A}_{\leq b}\subset\mathscr{A},
\end{equation*}
and we may intersect the subcategories of ingressive and egressive morphisms described to obtain the following exact $\infty$-category structure on $\mathscr{A}_{[a,b]}$. A morphism $\fromto{X}{Y}$ is ingressive just in case the induced morphism $\fromto{\pi_bX}{\pi_bY}$ is a monomorphism of the abelian category $\mathscr{A}^{\heartsuit}$. A morphism $\fromto{Y}{Z}$ is egressive just in case the induced morphism $\fromto{\pi_aX}{\pi_aY}$ is an epimorphism of $\mathscr{A}^{\heartsuit}$.
\end{enumerate}
\end{exm}

\begin{exm} Yet more generally, suppose $\mathscr{A}$ a stable $\infty$-category, and suppose $\mathscr{C}\subset\mathscr{A}$ \emph{any} full additive subcategory that is closed under extensions. Declare a morphism $\fromto{X}{Y}$ of $\mathscr{C}$ to be ingressive just in case its cofiber in $\mathscr{A}$ lies in $\mathscr{C}$. Dually, declare a morphism $\fromto{Y}{Z}$ of $\mathscr{C}$ to be egressive just in case its fiber in $\mathscr{A}$ lies in $\mathscr{C}$. Then $\mathscr{C}$ is exact with this triple structure.
\end{exm}

\begin{rmk} Thomason and Trobaugh called a triple $(C,C_{\dag},C^{\dag})$ of ordinary categories with direct sums whose nerves satisfy conditions (\ref{dfn:biWald}.\ref{item:exactinfty}.2-4) a \emph{category with bifibrations} \cite[Df. 1.2.2]{MR92f:19001}. Include additivity, and this is precisely the notion of an ordinary exact category. (As an aside, we remark that the theory of exact $\infty$-categories we delve into here really uses the additivity condition; without it, one would be unable to ensure that exact $\infty$-categories form a full subcategory of Waldhausen $\infty$-categories.)

When a class of weak equivalences is included, Thomason and Trobaugh used the term \emph{biWaldhausen category} \cite[Df. 1.2.4]{MR92f:19001}. This notion still does not require additivity. However, suppose $\mathbf{A}$ a \emph{complicial biWaldhausen category} \cite[Df. 1.2.11]{MR92f:19001} that is closed under the canonical homotopy pullbacks and canonical homotopy pushouts of \cite[Df. 1.1.2]{MR92f:19001} in the sense of \cite[Df. 1.3.5]{MR92f:19001} such that the mapping cylinder and cocylinder functors of \cite[Df. 1.3.4]{MR92f:19001} satisfy \cite[Df. 1.3.1.5]{MR92f:19001} and its dual. Then the relative nerve \cite[Df. 1.5]{K1} of $(\mathbf{A},w\mathbf{A})$ is clearly a stable $\infty$-category. Consequently, \emph{every} example of a category with cofibrations and weak equivalences that Thomason and Trobaugh study \cite[\S\S 2--11]{MR92f:19001} is actually nothing more than a model for some stable $\infty$-category.
\end{rmk}


\section{Exact functors between exact $\infty$-categories}

\begin{dfn} Suppose $\mathscr{C}$ and $\mathscr{D}$ two exact $\infty$-categories. A functor $F\colon\fromto{\mathscr{C}}{\mathscr{D}}$ will be said to be \emph{exact} if it preserves both cofibrations and fibrations and if $F$ is exact both as a functor of Waldhausen $\infty$-categories and as a functor of coWaldhausen $\infty$-categories.

We denote by $\Fun_{\Exact_{\infty}}(\mathscr{C},\mathscr{D})$ the full subcategory of $\Fun(\mathscr{C},\mathscr{D})$ spanned by the exact functors $\fromto{\mathscr{C}}{\mathscr{D}}$.
\end{dfn}
\noindent This definition, when set against with the definition of exact functor of Waldhausen categories (Rec. \ref{dfn:Waldinftycat}), appears to overburden the phrase ``exact functor'' and to create the possibility for some ambiguity; however, in Pr. \ref{lem:exactisexact} we will see that in fact no ambiguity obtains.

For now, let us construct the  $\infty$-category of exact $\infty$-categories.

\begin{ntn} Denote by $\Exact^{\Delta}_{\infty}$ the following simplicial category. The objects of $\Exact^{\Delta}_{\infty}$ are small exact $\infty$-categories; for any two exact $\infty$-categories $\mathscr{C}$ and $\mathscr{D}$, let $\Exact_{\infty}^{\Delta}(\mathscr{C},\mathscr{D})$ be the maximal Kan complex $\iota\Fun_{\Exact_{\infty}}(\mathscr{C},\mathscr{D})$ contained in $\Fun_{\Exact_{\infty}}(\mathscr{C},\mathscr{D})$. We write $\Exact_{\infty}$ for the simplicial nerve \cite[Df. 1.1.5.5]{HTT} of $\Exact_{\infty}^{\Delta}$.
\end{ntn}

\begin{rmk} We could have equally well defined $\Exact_{\infty}$ as the full subcategory of the pullback
\begin{equation*}
\Wald_{\infty}\times_{\Cat_{\infty}}\coWald_{\infty}
\end{equation*}
spanned by the exact $\infty$-categories.
\end{rmk}

We have already remarked that the formation of the opposite of a Waldhausen $\infty$-category defines an equivalence
\begin{equation*}
\equivto{\Wald_{\infty}}{\coWald_{\infty}}.
\end{equation*}
Since exact $\infty$-categories are defined by fitting together the structure of a Waldhausen $\infty$-category and a coWaldhausen $\infty$-category in a self-dual manner, we obtain the following.
\begin{lem}\label{lm:oppositeexact} The formation of the opposite restricts to an autoequivalence
\begin{equation*}
\op\colon\equivto{\Exact_{\infty}}{\Exact_{\infty}}.
\end{equation*}
\end{lem}
\noindent This permits us to dualize virtually any assertion about exact $\infty$-categories.

We now set about showing that the inclusion $\into{\Exact_{\infty}}{\Wald_{\infty}}$ is fully faithful. For this, we use in a nontrivial way the additivity condition for exact $\infty$-categories. In particular, this additivity actually guarantees a greater compatibility between pullbacks and pushouts and between fibrations and cofibrations than one might at first expect.
\begin{lem} In an exact $\infty$-category, a pushout square
\begin{equation*}
\begin{tikzpicture}[baseline]
\matrix(m)[matrix of math nodes, 
row sep=4ex, column sep=4ex, 
text height=1.5ex, text depth=0.25ex] 
{X&Y\\
X'&Y',\\}; 
\path[>=stealth,->,font=\scriptsize] 
(m-1-1) edge[>->] (m-1-2) 
edge (m-2-1) 
(m-1-2) edge (m-2-2) 
(m-2-1) edge[->] (m-2-2); 
\end{tikzpicture}
\end{equation*}
in which the morphism $\cofto{X}{Y}$ is ingressive is also a pullback square. Dually, a pullback square
\begin{equation*}
\begin{tikzpicture}[baseline]
\matrix(m)[matrix of math nodes, 
row sep=4ex, column sep=4ex, 
text height=1.5ex, text depth=0.25ex] 
{X&Y\\
X'&Y',\\}; 
\path[>=stealth,->,font=\scriptsize] 
(m-1-1) edge (m-1-2) 
edge (m-2-1) 
(m-1-2) edge[->>] (m-2-2) 
(m-2-1) edge (m-2-2); 
\end{tikzpicture}
\end{equation*}
in which the morphism $\fibto{Y}{Y'}$ is egressive is also a pushout square.
\begin{proof} We prove the first statement; the second is dual. Since $\cofto{X'}{Y'}$ is ingressive, we may form the cofiber
\begin{equation*}
\begin{tikzpicture}[baseline]
\matrix(m)[matrix of math nodes, 
row sep=4ex, column sep=4ex, 
text height=1.5ex, text depth=0.25ex] 
{X'&Y'\\
0&Z',\\}; 
\path[>=stealth,->,font=\scriptsize] 
(m-1-1) edge[>->] (m-1-2) 
edge[->>] (m-2-1) 
(m-1-2) edge[->>] (m-2-2) 
(m-2-1) edge[>->] (m-2-2); 
\end{tikzpicture}
\end{equation*}
which is an ambigressive square. Hence the square
\begin{equation*}
\begin{tikzpicture}[baseline]
\matrix(m)[matrix of math nodes, 
row sep=4ex, column sep=4ex, 
text height=1.5ex, text depth=0.25ex] 
{X&Y\\
0&Z',\\}; 
\path[>=stealth,->,font=\scriptsize] 
(m-1-1) edge[>->] (m-1-2) 
edge[->>] (m-2-1) 
(m-1-2) edge[->>] (m-2-2) 
(m-2-1) edge[>->] (m-2-2); 
\end{tikzpicture}
\end{equation*}
is also ambigressive, whence we conclude that 
\begin{equation*}
\begin{tikzpicture}[baseline]
\matrix(m)[matrix of math nodes, 
row sep=4ex, column sep=4ex, 
text height=1.5ex, text depth=0.25ex] 
{X&Y\\
X'&Y',\\}; 
\path[>=stealth,->,font=\scriptsize] 
(m-1-1) edge[>->] (m-1-2) 
edge (m-2-1) 
(m-1-2) edge (m-2-2) 
(m-2-1) edge[>->] (m-2-2); 
\end{tikzpicture}
\end{equation*}
is a pullback square.
\end{proof}
\end{lem}

The next pair of lemmas give a convenient way to replace pushout squares with exact sequences.

\begin{lem} For any exact sequence
\begin{equation*}
X'\ \tikz[baseline]\draw[>=stealth,>->,font=\scriptsize](0,0.5ex)--node[above]{$i$}(0.5,0.5ex);\ X\ \tikz[baseline]\draw[>=stealth,->>,font=\scriptsize](0,0.5ex)--node[above]{$p$}(0.5,0.5ex);\ X''
\end{equation*}
of an exact $\infty$-category $C$, the object $W$ formed as the pushout
\begin{equation*}
\begin{tikzpicture} 
\matrix(m)[matrix of math nodes, 
row sep=4ex, column sep=4ex, 
text height=1.5ex, text depth=0.25ex] 
{X'&X\\ 
X&W\\}; 
\path[>=stealth,>->,font=\scriptsize] 
(m-1-1) edge (m-1-2) 
edge (m-2-1) 
(m-1-2) edge (m-2-2) 
(m-2-1) edge (m-2-2); 
\end{tikzpicture}
\end{equation*}
is a direct sum $X\oplus X''$. Dually, the object $V$ formed as the pullback
\begin{equation*}
\begin{tikzpicture} 
\matrix(m)[matrix of math nodes, 
row sep=4ex, column sep=4ex, 
text height=1.5ex, text depth=0.25ex] 
{V&X\\ 
X&X''\\}; 
\path[>=stealth,->>,font=\scriptsize] 
(m-1-1) edge (m-1-2) 
edge (m-2-1) 
(m-1-2) edge (m-2-2) 
(m-2-1) edge (m-2-2); 
\end{tikzpicture}
\end{equation*}
is a direct sum $X'\oplus X$.
\begin{proof} We prove the first assertion; the second is dual. Choose a fibrant simplicial category $D$ whose nerve is equivalent to $C$. Now for any object $T$, the shear map
\begin{equation*}
\equivto{\Map_D(X,T)\times\Map_D(X,T)}{\Map_D(X,T)\times\Map_D(X,T)}
\end{equation*}
induces an equivalence
\begin{equation}\label{eq:dirsumaspush}
\begin{tikzpicture}[baseline]
\matrix(m)[matrix of math nodes, 
row sep=6ex, column sep=4ex, 
text height=1.5ex, text depth=0.25ex] 
{\left(\Map_D(X,T)\times\Map_D(X,T)\right)\times^{h}_{\id\times i^{\star},\left(\Map_D(X,T)\times\Map_D(X',T)\right),\id\times 0}\left(\Map_D(X,T)\times\Delta^0\right)\\ 
\left(\Map_D(X,T)\times\Map_D(X,T)\right)\times^{h}_{i^{\star}\times i^{\star},\left(\Map_D(X',T)\times\Map_D(X',T)\right),\Delta}\Map_D(X',T),\\}; 
\path[>=stealth,->,font=\scriptsize,inner sep=0.5pt] 
(m-1-1) edge node[right]{$\sim$} (m-2-1); 
\end{tikzpicture}
\end{equation}
where:
\begin{itemize}
\item[---] $i^{\star}$ denotes the map $\fromto{\Map_D(X,T)}{\Map_D(Z,T)}$ induced by $i\colon\cofto{X'}{X}$,
\item[---] $0$ denotes a vertex $\fromto{\Delta^0}{\Map_D(X',T)}$ corresponding to a zero map, and
\item[---] $\Delta$ denotes the diagonal map.
\end{itemize}
The source of \eqref{eq:dirsumaspush} is the product of $\Map_D(X,T)$ with the space of squares of the form
\begin{equation*}
\begin{tikzpicture} 
\matrix(m)[matrix of math nodes, 
row sep=4ex, column sep=4ex, 
text height=1.5ex, text depth=0.25ex] 
{X'&X\\ 
0&T,\\}; 
\path[>=stealth,->,font=\scriptsize] 
(m-1-1) edge[>->] node[above]{$i$} (m-1-2) 
edge (m-2-1) 
(m-1-2) edge (m-2-2) 
(m-2-1) edge (m-2-2); 
\end{tikzpicture}
\end{equation*}
in $C$, and the target is equivalent to the space of squares of the form
\begin{equation*}
\begin{tikzpicture} 
\matrix(m)[matrix of math nodes, 
row sep=4ex, column sep=4ex, 
text height=1.5ex, text depth=0.25ex] 
{X'&X\\ 
X&T\\}; 
\path[>=stealth,->,font=\scriptsize] 
(m-1-1) edge[>->] node[above]{$i$} (m-1-2) 
edge[>->] node[left]{$i$} (m-2-1) 
(m-1-2) edge (m-2-2) 
(m-2-1) edge (m-2-2); 
\end{tikzpicture}
\end{equation*}
in $C$. Consequently, the map \eqref{eq:dirsumaspush} specifies an equivalence
\begin{equation*}
\equivto{\Map_{D}(X,T)\times\Map_{D}(X'',T)}{\Map(W,T)}.
\end{equation*}
This equivalence is clearly functorial in $T$, so it specifies an equivalence $\equivto{W}{X\oplus X''}$.
\end{proof}
\end{lem}

\begin{lem} In an exact $\infty$-category, suppose that
\begin{equation*}
\begin{tikzpicture}[baseline]
\matrix(m)[matrix of math nodes, 
row sep=4ex, column sep=4ex, 
text height=1.5ex, text depth=0.25ex] 
{X&Y\\
X'&Y',\\}; 
\path[>=stealth,->,font=\scriptsize] 
(m-1-1) edge node[above]{$i$} (m-1-2) 
edge node[left]{$p$} (m-2-1) 
(m-1-2) edge node[right]{$q$} (m-2-2) 
(m-2-1) edge node[below]{$i'$} (m-2-2); 
\end{tikzpicture}
\end{equation*}
is either a pushout square in which $\cofto{X}{Y}$ is ingressive or a pullback square in which $\fibto{Y}{Y'}$ is egressive. Then the morphism
\begin{equation*}
{-p\choose i}\colon\cofto{X}{X'\oplus Y}
\end{equation*}
is ingressive, the morphism
\begin{equation*}
(i'\quad q)\colon\fibto{X'\oplus Y}{Y'}
\end{equation*}
is egressive, and these maps fit into an exact sequence
\begin{equation*}
X\ \tikz[baseline]\draw[>=stealth,>->](0,0.5ex)--(0.5,0.5ex);\ X'\oplus Y\ \tikz[baseline]\draw[>=stealth,->>](0,0.5ex)--(0.5,0.5ex);\ Y'.
\end{equation*}
\begin{proof} We prove the assertion for pushout squares; the other assertion is dual. We form a diagram
\begin{equation*}
\begin{tikzpicture} 
\matrix(m)[matrix of math nodes, 
row sep=4ex, column sep=4ex, 
text height=1.5ex, text depth=0.25ex] 
{X&X'&0\\ 
V&V'&Y'\\
Y&Y'&Z\\}; 
\path[>=stealth,->,font=\scriptsize] 
(m-1-1) edge (m-1-2) 
edge[>->] (m-2-1)
(m-1-2) edge[->>] (m-1-3)
edge[>->] (m-2-2)
(m-1-3) edge[>->] (m-2-3)
(m-2-1) edge (m-2-2)
edge[->>] (m-3-1)
(m-2-2) edge[->>] (m-2-3)
edge[->>] (m-3-2)
(m-2-3) edge[->>] (m-3-3)
(m-3-1) edge (m-3-2)
(m-3-2) edge[->>] (m-3-3); 
\end{tikzpicture}
\end{equation*}
in which every square is a pullback square. By the previous lemma, $V'$ is a direct sum $X'\oplus Y'$, and $V$ is a direct sum $X\oplus Y'$. The desired exact sequence is the top rectangle.
\end{proof}
\end{lem}

\begin{prp}\label{lem:exactisexact} The following are equivalent for a functor $\psi\colon\fromto{\mathscr{C}}{\mathscr{D}}$ between two $\infty$-categories with exact $\infty$-category structures.
\begin{enumerate}[(\ref{lem:exactisexact}.1)]
\item The functor $\psi$ carries cofibrations to cofibrations, it carries fibrations to fibrations, and as a functor of exact $\infty$-categories, $\psi$ is exact.
\item The functor $\psi$ carries cofibrations to cofibrations, and as a functor of Waldhausen $\infty$-categories, $\psi$ is exact.
\item The functor $\psi$ carries fibrations to fibrations, and as a functor of coWaldhausen $\infty$-categories, $\psi$ is exact.
\end{enumerate}
\begin{proof} It is clear that the first condition implies the other two. We shall show that the second implies the first; the proof that the third condition implies the first is dual. So suppose $\psi$ preserves cofibrations and is exact as a functor of Waldhausen $\infty$-categories. Because a morphism is egressive just in case it can be exhibited as a cofiber, $\psi$ preserves fibrations as well. A pullback square
\begin{equation*}
\begin{tikzpicture} 
\matrix(m)[matrix of math nodes, 
row sep=4ex, column sep=4ex, 
text height=1.5ex, text depth=0.25ex] 
{X&Y\\ 
X'&Y'\\}; 
\path[>=stealth,->,font=\scriptsize] 
(m-1-1) edge node[above]{$i$} (m-1-2) 
edge[->>] node[left]{$p$} (m-2-1) 
(m-1-2) edge[->>] node[right]{$q$} (m-2-2) 
(m-2-1) edge node[below]{$i'$} (m-2-2); 
\end{tikzpicture}
\end{equation*}
in which $p\colon\fibto{X}{X'}$ and $q\colon\fibto{Y}{Y'}$ are egressive can be factored as
\begin{equation*}
\begin{tikzpicture} 
\matrix(m)[matrix of math nodes, 
row sep=6ex, column sep=6ex, 
text height=1.5ex, text depth=0.25ex] 
{X&X'\oplus Y&Y\\ 
X'&X'\oplus Y'&Y'\\}; 
\path[>=stealth,->,font=\scriptsize] 
(m-1-1) edge[>->] node[above]{$-p\choose i$} (m-1-2) 
edge[->>] node[left]{$p$} (m-2-1) 
(m-1-2) edge[->>] node[right]{$\id\oplus q$} (m-2-2)
edge[->>] node[above]{$\pr_2$} (m-1-3)
(m-1-3) edge[->>] node[right]{$q$} (m-2-3)
(m-2-1) edge[>->] node[below]{$-\id\choose i'$} (m-2-2)
(m-2-2) edge[->>] node[below]{$\pr_2$} (m-2-3); 
\end{tikzpicture}
\end{equation*}
in which all three rectangles are pullbacks. By the previous lemma, the left hand square is an ambigressive pullback/pushout, so when we apply $\psi$, we obtain a diagram
\begin{equation*}
\begin{tikzpicture} 
\matrix(m)[matrix of math nodes, 
row sep=6ex, column sep=6ex, 
text height=1.5ex, text depth=0.25ex] 
{\psi X&\psi X'\oplus\psi Y&\psi Y\\ 
\psi X'&\psi X'\oplus\psi Y'&\psi Y'\\}; 
\path[>=stealth,->,font=\scriptsize] 
(m-1-1) edge[>->] node[above]{$-\psi p\choose \psi i$} (m-1-2) 
edge[->>] node[left]{$\psi p$} (m-2-1) 
(m-1-2) edge[->>] node[right]{$\id\oplus \psi q$} (m-2-2)
edge[->>] node[above]{$\pr_2$} (m-1-3)
(m-1-3) edge[->>] node[right]{$\psi q$} (m-2-3)
(m-2-1) edge[>->] node[below]{$-\id\choose \psi i'$} (m-2-2)
(m-2-2) edge[->>] node[below]{$\pr_2$} (m-2-3); 
\end{tikzpicture}
\end{equation*}
in which the right hand square is easily seen to be a pullback, and the left hand square, being an ambigressive pushout, is also an ambigressive pullback.
\end{proof}
\end{prp}
\begin{cor}\label{cor:fullfaith} The forgetful functors
\begin{equation*}
\fromto{\Exact_{\infty}}{\Wald_{\infty}}\textrm{\quad and\quad}\fromto{\Exact_{\infty}}{\coWald_{\infty}}
\end{equation*}
are fully faithful.
\end{cor}
\noindent In particular, we may say that a Waldhausen $\infty$-category $\mathscr{C}$ ``is'' an exact $\infty$-category if it lies in the essential image of the forgetful functor $\fromto{\Exact_{\infty}}{\Wald_{\infty}}$, and we will treat this forgetful functor as if it were an inclusion. Since this functor is fully faithful, this is not a significant abuse of terminology. We make sense of the assertion that a coWaldhausen $\infty$-category ``is'' an exact $\infty$-category in a dual manner.

\begin{nul}\label{nul:essimofforExacttoWald} The essential image of the forgetful functor $\fromto{\Exact_{\infty}}{\Wald_{\infty}}$ is spanned by those Waldhausen $\infty$-categories that satisfy the following three criteria.
\begin{enumerate}[(\ref{nul:essimofforExacttoWald}.1)]
\item The underlying $\infty$-category is additive.
\item The class of morphisms that can be exhibited as the cofiber of some cofibration is closed under pullback.
\item Every cofibration is the fiber of its cofiber.
\end{enumerate}
The essential image of the forgetful functor $\fromto{\Exact_{\infty}}{\coWald_{\infty}}$ is described in a dual manner.
\end{nul}

\section{Theories and duality}\label{sect:duality} Algebraic $K$-theory is a particular example of what we called an additive  theory in \cite{K1}. In effect, additive theories are the natural homology theories for Waldhausen $\infty$-categories. To tell this story, it is necessary to recall some pleasant facts about the $\infty$-categories $\Wald_{\infty}$ and $\coWald_{\infty}$.

\begin{rec} In particular, the $\infty$-category $\Wald_{\infty}$ (and hence also the $\infty$-category $\coWald_{\infty}$) enjoys a number of excellent formal properties. We showed in \cite[Pr. 4.6]{K1} that it admits direct sums, and we also showed in \cite[Pr. 4.7]{K1} that it is compactly generated \cite[Df. 5.5.7.1]{HTT} in the sense that every Waldhausen $\infty$-category is in fact the filtered union of its finitely presented Waldhausen subcategories (that is, of Waldhausen subcategories that are compact as objects of $\Wald_{\infty}$).

Furthermore, suppose $\mathscr{X}\colon\fromto{\Lambda}{\Wald_{\infty}}$ a diagram of Waldhausen $\infty$-categories. The limit $\lim\mathscr{X}$ is computed by forming the limit in $\Cat_{\infty}$ \cite[Cor. 3.3.3.2]{HTT} and then declaring a morphism to be ingressive if its image in each $\mathscr{X}_{\alpha}$ is so \cite[Pr. 4.3]{K1}. Similarly, if $\Lambda$ is filtered, then the colimit $\colim\mathscr{X}$ is computed by forming the colimit in $\Cat_{\infty}$ \cite[Cor. 3.3.4.3]{HTT} and then declaring $(\colim\mathscr{X})_{\dag}$ to be the union of the images of the $\infty$-categories $\mathscr{X}_{\alpha,\dag}$ \cite[Pr. 4.4]{K1}.

Since the $\infty$-categories $\Wald_{\infty}$ and $\coWald_{\infty}$ are equivalent, it is clear that all these properties are enjoyed by the latter as well as the former. We'll denote by $\Wald_{\infty}^{\omega}$ (respectively, by $\coWald_{\infty}^{\omega}$) the full subcategory of $\Wald_{\infty}$ (resp., $\coWald_{\infty}$) spanned by the finitely presented Waldhausen $\infty$-categories (resp., the finitely presented coWaldhausen $\infty$-categories).
\end{rec}

These formal properties can be regarded as analogues of a few of the formal properties enjoyed by the ordinary category $\VV(k)$ of vector spaces: vector spaces are of course additive, and any vector space is the union of its finite-dimensional subspaces. Furthermore, the underlying set of a limit or filtered colimit of vector spaces is the limit or filtered colimit of the underlying sets.

We will only be interested in functors on $\Wald_{\infty}$ or $\coWald_{\infty}$ that are (1) trivial on the zero Waldhausen $\infty$-category, and (2) are completely determined by their values on finitely presented  Waldhausen $\infty$-categories. 

\begin{rec}[\protect{\cite[Df. 7.1]{K1}}] A \emph{theory} $\phi\colon\fromto{\Wald_{\infty}}{\Kan}$ is a functor that preserves terminal objects and filtered colimits.

Similarly, a \emph{theory} $\phi\colon\fromto{\coWald_{\infty}}{\Kan}$ is a functor that preserves terminal objects and filtered colimits.
\end{rec}

\begin{ntn} Denote by $\Thy$ (respectively, $\Thy^{\vee}$) the full subcategory of the $\infty$-category $\Fun(\coWald_{\infty},\Kan)$ spanned by the theories.
\end{ntn}

\begin{exm} The most important example of a theory (in either sense) is the \emph{moduli space of objects} functor $\goesto{\mathscr{C}}{\iota\mathscr{C}}$. Here $\iota\mathscr{C}$ denotes the largest Kan complex contained in $\mathscr{C}$ \cite[Pr. 1.2.5.3]{HTT}.
\end{exm}

Of course there is a canonical equivalence between theories on Waldhausen $\infty$-categories and theories on coWaldhausen $\infty$-categories.

\begin{dfn} For any theory $\phi\colon\fromto{\Wald_{\infty}}{\Kan}$, the \emph{dual theory} $\phi^{\vee}$ is the composite $\phi\circ\op\colon\fromto{\coWald_{\infty}}{\Kan}$. This construction clearly yields an equivalence of $\infty$-categories $\equivto{\Thy}{\Thy^{\vee}}$.
\end{dfn}

Now we are prepared to define the key notion of an \emph{exact duality} on a theory $\phi$.

\begin{dfn}\label{dfn:theoryforexacts} An \emph{exact duality} on a theory $\phi$ is an equivalence
\begin{equation*}
\eta\colon\equivto{\phi|_{\Exact_{\infty}}}{\phi^{\vee}|_{\Exact_{\infty}}}
\end{equation*}
of the $\infty$-category $\Fun(\Exact_{\infty},\mathscr{E})$.
\end{dfn}

\begin{exm} Suppose $\rho\colon\fromto{\Cat_{\infty}^{\ast}}{\Kan}$ a functor from the $\infty$-category $\Cat_{\infty}^{\ast}$ of pointed $\infty$-categories to $\Kan$ that preserves the terminal object and filtered colimits. Then an equivalence $\equivto{\rho}{\rho\circ\op}$ induces an exact duality on the composite
\begin{equation*}
\Wald_{\infty}\to\Cat_{\infty}^{\ast}\to\Kan.
\end{equation*}
For instance, the functor $\iota\colon\fromto{\Cat_{\infty}^{\ast}}{\Kan}$ admits an equivalence $\equivto{\iota}{\iota\circ\op}$; consequently, the theory $\iota\colon\fromto{\Wald_{\infty}}{\Kan}$ admits an exact duality.
\end{exm}

A general theory doesn't reflect much about the (co)Waldhausen structure. \emph{Additive} theories are much more sensitive. To talk about them, we have to recall our construction of $\mathscr{F}$ and $\mathscr{S}$.
\begin{rec} In \cite[\S 5]{K1} we defined, for any Waldhausen $\infty$-category $\mathscr{C}$, an $\infty$-category $\mathscr{F}(\mathscr{C})$ and a full subcategory $\mathscr{S}(\mathscr{C})\subset\mathscr{F}(\mathscr{C})$. An object of $\mathscr{F}(\mathscr{C})$ is a pair $(m,X)$ consisting of a natural number $m$ and a sequence
\begin{equation*}
X_0\ \tikz[baseline]\draw[>=stealth,>->](0,0.5ex)--(0.5,0.5ex);\ \cdots\ \tikz[baseline]\draw[>=stealth,>->](0,0.5ex)--(0.5,0.5ex);\ X_m
\end{equation*}
of cofibrations of $\mathscr{C}$; a morphism $\fromto{(m,X)}{(n,Y)}$ of $\mathscr{F}(\mathscr{C})$ consists of a morphism $\eta\colon\fromto{\mathbf{n}}{\mathbf{m}}$ of $\Delta$ and a commutative diagram
\begin{equation*}
\begin{tikzpicture} 
\matrix(m)[matrix of math nodes, 
row sep=4ex, column sep=4ex, 
text height=1.5ex, text depth=0.25ex] 
{X_{\eta(0)}&\cdots&X_{\eta(n)}\\ 
Y_0&\cdots&Y_n\\}; 
\path[>=stealth,->,font=\scriptsize] 
(m-1-1) edge[>->] (m-1-2) 
edge (m-2-1) 
(m-1-2) edge[>->] (m-1-3)
(m-1-3) edge (m-2-3)
(m-2-1) edge[>->] (m-2-2)
(m-2-2) edge[>->] (m-2-3); 
\end{tikzpicture}
\end{equation*}
of $\mathscr{C}$. The full subcategory $\mathscr{S}(\mathscr{C})\subset\mathscr{F}(\mathscr{C})$ is the one spanned by those pairs $(m,X)$ such that $X_0$ is a zero object of $\mathscr{C}$.

The assignment $\goesto{(m,X)}{\mathbf{m}}$ defines functors
\begin{equation*}
p_{\mathscr{C}}\colon\fromto{\mathscr{F}(\mathscr{C})}{N\Delta^{\op}}\quad and\quad q_{\mathscr{C}}\colon\fromto{\mathscr{S}(\mathscr{C})}{N\Delta^{\op}}.
\end{equation*}
Let us write $\mathscr{F}_m(\mathscr{C})\coloneq p_{\mathscr{C}}^{-1}(\mathbf{m})$ and $\mathscr{S}_m(\mathscr{C})\coloneq q_{\mathscr{C}}^{-1}(\mathbf{m})$.

The first basic fact to understand about these functors is that they are \emph{cocartesian fibrations} \cite[Df. 2.4.2.1]{HTT}; that is, there are functors
\begin{equation*}
\FF_{\ast}(\mathscr{C})\colon\fromto{N\Delta^{\op}}{\Cat_{\infty}}\textrm{\quad and\quad}\SS_{\ast}(\mathscr{C})\colon\fromto{N\Delta^{\op}}{\Cat_{\infty}}
\end{equation*}
(the functors that classify $p_{\mathscr{C}}$ and $q_{\mathscr{C}}$, \cite[Df. 3.3.2.2]{HTT}) and equivalences $\FF_m(\mathscr{C})\simeq\mathscr{F}_m(\mathscr{C})$ and $\SS_m(\mathscr{C})\simeq\mathscr{S}_m(\mathscr{C})$ such that for any morphism $\eta\colon\fromto{\mathbf{n}}{\mathbf{m}}$ of $\Delta$, the space of morphisms $\fromto{(m,X)}{(n,Y)}$ of $\mathscr{F}(\mathscr{C})$ (respectively, of $\mathscr{S}(\mathscr{C})$) that cover $\eta$ is equivalent to the space of morphisms $\fromto{\eta_!X}{Y}$, where $\eta_!$ is shorthand notation for the image of $\eta$ under $\FF_{\ast}(\mathscr{C})$ (resp., $\SS_{\ast}(\mathscr{C})$). In particular, for any object $X$ of $p_{\mathscr{C}}^{-1}(\mathbf{m})$ (resp., $q_{\mathscr{C}}^{-1}(\mathbf{m})$), there exists a special edge --- called a \emph{cocartesian edge} --- from $X$ to $\eta_!X$.

The functor $\FF_{\ast}(\mathscr{C})$ is easy to describe: it carries $\mathbf{m}$ to the full subcategory of $\Fun(\Delta^m,\mathscr{C})$ spanned by those functors $X\colon\fromto{\Delta^m}{\mathscr{C}}$ such that each morphism $\fromto{X_i}{X_{i+1}}$ is ingressive; the functoriality in $\mathbf{m}$ is obvious here. The functor $\SS_{\ast}(\mathscr{C})$ is a tad trickier to describe: morally, it carries $\mathbf{m}$ to the full subcategory of $\FF_m(\mathscr{C})$ spanned by those objects $X$ such that $X_0$ is a zero object, and a map $\eta\colon\fromto{\mathbf{n}}{\mathbf{m}}$ of $\Delta$ induces a functor $\fromto{\SS_m(\mathscr{C})}{\SS_n(\mathscr{C})}$ that carries an object $X$ to the object
\begin{equation*}
0\ \tikz[baseline]\draw[>=stealth,>->](0,0.5ex)--(0.5,0.5ex);\ X_{\eta(1)}/X_{\eta(0)}\ \tikz[baseline]\draw[>=stealth,>->](0,0.5ex)--(0.5,0.5ex);\ \cdots\ \tikz[baseline]\draw[>=stealth,>->](0,0.5ex)--(0.5,0.5ex);\ X_{\eta(n)}/X_{\eta(0)}.
\end{equation*}
One can opt to make compatible choices of these quotients to rectify this into an actual functor of $\infty$-categories, or, alternately, one can forget about $\SS_{\ast}(\mathscr{C})$ and stick with the cocartesian fibration $\fromto{\mathscr{S}(\mathscr{C})}{N\Delta^{\op}}$. In \cite{MR86m:18011}, Waldhausen opted for the former; in \cite{K1}, we opted for the latter. Which approach one chooses to take is largely a matter of taste or convenience; conceptually, there is no difference.

In \cite[Cor. 5.20.1]{K1}, we show that the left adjoints $\fromto{\mathscr{F}_m(\mathscr{C})}{\mathscr{S}_m(\mathscr{C})}$ to the inclusions $\into{\mathscr{S}_m(\mathscr{C})}{\mathscr{F}_m\mathscr{C}}$ fit together over $N\Delta^{\op}$ to make a natural transformation $\fromto{\FF_{\ast}(\mathscr{C})}{\SS_{\ast}(\mathscr{C})}$.

The next significant thing to understand about $\fromto{\mathscr{F}(\mathscr{C})}{N\Delta^{\op}}$ and $\fromto{\mathscr{S}(\mathscr{C})}{N\Delta^{\op}}$ is that they are what we have called \emph{Waldhausen cocartesian fibrations} \cite[Df. 3.20]{K1}. That is, there are subcategories
\begin{equation*}
\mathscr{F}(\mathscr{C})_{\dag}\subset\mathscr{F}(\mathscr{C})\textrm{\quad and\quad}\mathscr{S}(\mathscr{C})_{\dag}\subset\mathscr{S}(\mathscr{C})
\end{equation*}
such that for any $m\geq 0$, the pairs
\begin{equation*}
(\mathscr{F}_m(\mathscr{C}),\mathscr{F}_m(\mathscr{C})\cap\mathscr{F}(\mathscr{C})_{\dag})\textrm{\quad and\quad}(\mathscr{S}_m(\mathscr{C}),\mathscr{S}_m(\mathscr{C})\cap\mathscr{S}(\mathscr{C})_{\dag})
\end{equation*}
are Waldhausen $\infty$-categories, and for any morphism $\eta\colon\fromto{\mathbf{n}}{\mathbf{m}}$ of $\Delta$, the functors
\begin{equation*}
\eta_!\colon\fromto{\mathscr{F}_m(\mathscr{C})\simeq\FF_m(\mathscr{C})}{\FF_n(\mathscr{C})\simeq\mathscr{F}_n(\mathscr{C})}\textrm{\quad and\quad}\eta_!\colon\fromto{\mathscr{S}_m(\mathscr{C})\simeq\SS_m(\mathscr{C})}{\SS_n(\mathscr{C})\simeq\mathscr{S}_n(\mathscr{C})}
\end{equation*}
are exact functors of Waldhausen $\infty$-categories. In effect, a morphism $\fromto{(m,X)}{(n,Y)}$ of either $\mathscr{F}(\mathscr{C})$ or $\mathscr{S}(\mathscr{C})$ will be declared ingressive just in case the morphism $\fromto{\mathbf{n}}{\mathbf{m}}$ is an isomorphism and the diagram
\begin{equation*}
\begin{tikzpicture} 
\matrix(m)[matrix of math nodes, 
row sep=4ex, column sep=4ex, 
text height=1.5ex, text depth=0.25ex] 
{X_0&\cdots&X_m\\ 
Y_0&\cdots&Y_m\\}; 
\path[>=stealth,->,font=\scriptsize] 
(m-1-1) edge[>->] (m-1-2) 
edge (m-2-1) 
(m-1-2) edge[>->] (m-1-3)
(m-1-3) edge (m-2-3)
(m-2-1) edge[>->] (m-2-2)
(m-2-2) edge[>->] (m-2-3); 
\end{tikzpicture}
\end{equation*}
has the property that for any integer $0\leq i\leq j\leq m$, the natural morphism $\fromto{Y_i\cup^{X_i}X_j}{Y_j}$ (whose source is the pushout in the $\infty$-categorical sense, of course) is ingressive. We show that $\fromto{\mathscr{F}(\mathscr{C})}{N\Delta^{\op}}$ and $\fromto{\mathscr{S}(\mathscr{C})}{N\Delta^{\op}}$ are Waldhausen cocartesian fibrations in \cite[Pr. 5.11 and Th. 5.20]{K1}. Consequently, the functors that classify these fibrations can be lifted to simplicial Waldhausen $\infty$-categories
\begin{equation*}
\FF_{\ast}(\mathscr{C})\colon\fromto{N\Delta^{\op}}{\Wald_{\infty}}\textrm{\quad and\quad}\SS_{\ast}(\mathscr{C})\colon\fromto{N\Delta^{\op}}{\Wald_{\infty}}
\end{equation*}

Dually, for any coWaldhausen $\infty$-category $\mathscr{D}$, we obtain \emph{coWaldhausen cartesian fibrations}
\begin{equation*}
\fromto{\mathscr{F}^{\vee}(\mathscr{D})\coloneq\mathscr{F}(\mathscr{D}^{\op})^{\op}}{N\Delta}\textrm{\quad and\quad}\fromto{\mathscr{S}^{\vee}(\mathscr{D})\coloneq\mathscr{S}(\mathscr{D}^{\op})^{\op}}{N\Delta}.
\end{equation*}
The functors that classify these fibrations are simplicial coWaldhausen $\infty$-categories
\begin{equation*}
\FF^{\vee}_{\ast}(\mathscr{D})\colon\fromto{N\Delta^{\op}}{\coWald_{\infty}}\textrm{\quad and\quad}\SS^{\vee}_{\ast}(\mathscr{D})\colon\fromto{N\Delta^{\op}}{\coWald_{\infty}}.
\end{equation*}
\end{rec}

\begin{dfn}[\protect{\cite[Df. 7.5]{K1}}] A theory $\phi\colon\fromto{\Wald_{\infty}}{\Kan}$ will be said to be \emph{additive} if for any Waldhausen $\infty$-category $\mathscr{C}$, the simplicial space
\begin{equation*}
\phi\circ\SS_{\ast}(\mathscr{C})\colon\fromto{N\Delta^{\op}}{\Kan}
\end{equation*}
is a group object in the sense of \cite[Df. 7.2.2.1]{HTT}. That is, for any $m\geq 0$, the map
\begin{equation*}
\fromto{\phi\circ\SS_{m}(\mathscr{C})}{\prod_{i=1}^m\phi\circ\SS_{\{i-1,i\}}(\mathscr{C})}
\end{equation*}
is an equivalence, and the monoid $\pi_0(\phi\circ\SS_1(\mathscr{C}))$ is a group. Write $\Add$ for the full subcategory of $\Thy$ spanned by the additive theories.

Dually, a theory $\phi\colon\fromto{\coWald_{\infty}}{\Kan}$ will be said to be \emph{additive} just in case it is the dual theory of an additive theory $\fromto{\Wald_{\infty}}{\Kan}$. That is, $\phi$ is additive just in case, for any coWaldhausen $\infty$-category $\mathscr{C}$, the simplicial space
\begin{equation*}
\phi\circ\SS^{\vee}_{\ast}(\mathscr{C})\colon\fromto{N\Delta^{\op}}{\Kan}
\end{equation*}
is a group object. Write $\Add^{\vee}$ for the full subcategory of $\Thy^{\vee}$ spanned by the additive theories.
\end{dfn}

\begin{rec} The main theorems of \cite{K1} are (1) that there exist what we call the \emph{fissile derived $\infty$-category} $\D_{\mathrm{add}}(\Wald_{\infty})$ and a functor
\begin{equation*}
\fromto{\Wald_{\infty}}{\D_{\mathrm{add}}(\Wald_{\infty})}
\end{equation*}
an equivalence
\begin{equation*}
\Exc(\D_{\mathrm{add}}(\Wald_{\infty}),\Kan)\simeq\Add,
\end{equation*}
where $\Exc$ denotes the $\infty$-category of reduced excisive functors that preserve sifted colimits \cite[Th. 7.4 and 7.6]{K1}, and (2) that the suspension in $\D_{\mathrm{add}}(\Wald_{\infty})$ of a Waldhausen $\infty$-category $\mathscr{C}$ is given by a formal colimit $\colim\SS_{\ast}(\mathscr{C})$ \cite[Cor. 6.9.1]{K1}. As a consequence, we deduced one may use the Goodwillie differential to find the best additive approximation to a theory $\phi$:
\begin{equation*}
D(\phi)(\mathscr{C})\simeq\underset{m\rightarrow\infty}{\colim}\ \Omega^m|\phi\circ\SS_{\ast}\cdots\SS_{\ast}(\mathscr{C})|_{(N\Delta^{\op})^m},
\end{equation*}
where $|\cdot|_{(N\Delta^{\op})^m}$ denotes the colimit of an $m$-simplicial space \cite[Th. 7.8]{K1}. The functor $D$ then defines a left adjoint to the inclusion $\into{\Thy}{\Add}$.

Dually, the inclusion $\into{\Add^{\vee}}{\Thy^{\vee}}$ admits a left adjoint $D^{\vee}$, given by the formula
\begin{equation*}
D^{\vee}(\phi)(\mathscr{C})\simeq\underset{m\rightarrow\infty}{\colim}\ \Omega^m|\phi\circ\SS_{\ast}^{\vee}\cdots\SS_{\ast}^{\vee}(\mathscr{C})|_{(N\Delta^{\op})^m}.
\end{equation*}
\end{rec}

The purpose of the remainder of this section is to describe a circumstance in which an exact duality on a theory $\phi$ descends to an exact duality on $D\phi$, and to show that these conditions obtain when $\phi=\iota$, giving a functorial equivalence $K(\mathscr{C})\simeq K(\mathscr{C}^{\op})$ for exact $\infty$-categories $\mathscr{C}$.

So suppose $\phi$ a theory with an exact duality. Note that $(D(\phi))^{\vee}$ is by construction equivalent to $D^{\vee}(\phi^{\vee})$, so such a result can be thought of as giving an equivalence $D(\phi)\simeq D^{\vee}(\phi^{\vee})$ on exact $\infty$-categories.

Consequently, for an exact $\infty$-category $\mathscr{C}$, we aim to produce a kind of duality between the Waldhausen cocartesian fibration $\fromto{\mathscr{S}(\mathscr{C})}{N\Delta^{\op}}$ and the coWaldhausen cartesian fibration $\fromto{\mathscr{S}^{\vee}(\mathscr{C})}{N\Delta}$.

More precisely, we will construct a functor $\widetilde{\SS}_{\ast}(\mathscr{C})\colon\fromto{N\Delta^{\op}}{\Cat_{\infty}}$ such that on the one hand, $\widetilde{\SS}_{\ast}$ classifies the cocartesian fibration $\fromto{\mathscr{S}(\mathscr{C})}{N\Delta^{\op}}$, and on the other, the composite of $\widetilde{\SS}_{\ast}$ with the functor $\op\colon\fromto{N\Delta^{\op}}{N\Delta^{\op}}$ given by $\goesto{\mathbf{[n]}}{\mathbf{[n]}^{\op}}$ is a straightening of the cartesian fibration $\fromto{\mathscr{S}^{\vee}(\mathscr{C})}{N\Delta}$.

In order to do this, we introduce thickened versions $\widetilde{\mathscr{S}},\ \widetilde{\mathscr{S}}^{\vee}$ of the constructions $\mathscr{S},\ \mathscr{S}^{\vee}$.

\begin{ntn} Let $\widetilde{\mathrm{M}}$ be the following ordinary category. The objects are triples $(m,i,j)$ consisting of integers $0\leq i\leq j\leq m$, and a morphism $\fromto{(n,k,\ell)}{(m,i,j)}$ is be a morphism $\phi\colon\fromto{\mathbf{[m]}}{\mathbf{[n]}}$ of $\Delta$ such that $k\leq\phi(i)$ and $\ell\leq\phi(j)$.

We declare an edge $\fromto{(n,k,\ell)}{(m,i,j)}$ of $N\widetilde{\mathrm{M}}$ to be \emph{ingressive} if the underlying edge $\fromto{\mathbf{m}}{\mathbf{n}}$ of $\Delta$ is an isomorphism and if $\ell=j$. Dually, we declare an edge $\fromto{(n,k,\ell)}{(m,i,j)}$ of  $N\widetilde{\mathrm{M}}$ to be egressive if the underlying edge $\fromto{\mathbf{m}}{\mathbf{n}}$ of $\Delta$ is an isomorphism and if $i=k$. We write $N\widetilde{\mathrm{M}}_{\dag}$ for the subcategory of $N\widetilde{\mathrm{M}}$ consisting of the ingressive morphisms, and we write $N\widetilde{\mathrm{M}}^{\dag}$ for the subcategory of $N\widetilde{\mathrm{M}}$ consisting of the egressive morphisms.

The fiber of the functor $\fromto{N\widetilde{\mathrm{M}}}{N\Delta^{\op}}$ over a vertex $\mathbf{n}\in N\Delta^{\op}$ is the arrow $\infty$-category $\mathscr{O}(\Delta^n)\coloneq\Fun(\Delta^1,\Delta^n)$.

Now the functor $\fromto{N\widetilde{\mathrm{M}}}{N\Delta^{\op}}$ is a cartesian fibration, and so its opposite $\fromto{N\widetilde{\mathrm{M}}^{\op}}{N\Delta}$ is a cocartesian fibration. Furthermore, one verifies easily that for any morphism $\eta\colon\fromto{\mathbf{n}}{\mathbf{m}}$, the induced functor
\begin{equation*}
\eta_{!}\colon\fromto{\mathscr{O}(\Delta^n)}{\mathscr{O}(\Delta^m)}
\end{equation*}
is the obvious one, and it preserves both ingressive and egressive morphisms.
\end{ntn}

\begin{cnstr} If $\mathscr{C}$ is a Waldhausen $\infty$-category, write $\widetilde{\mathscr{F}}(\mathscr{C})$ for the simplicial set over $N\Delta^{\op}$ satisfying the following universal property, which follows the general pattern set in \cite[Cor. 3.2.2.13]{HTT}. We require, for any simplicial set $K$ and any map $\sigma\colon\fromto{K}{N\Delta^{\op}}$, a bijection
\begin{equation*}
\Mor_{/(N\Delta^{\op})}(K,\widetilde{\mathscr{F}}(\mathscr{C}))\cong\Mor_{s\Set(2)}((K\times_{N\Delta^{\op}}N\widetilde{\mathrm{M}},K\times_{N\Delta^{\op}}(N\widetilde{\mathrm{M}})_{\dag}),(\mathscr{C},\mathscr{C}_{\dag})),
\end{equation*}
functorial in $\sigma$. Here, the category $s\Set(2)$ is the ordinary category of pairs $(X,A)$ of simplicial sets $X$ equipped with a simplicial subset $A\subset X$. (Note that since the functor of $K$ on the right hand side carries colimits to limits, the simplicial set $\widetilde{\mathscr{F}}(\mathscr{C})$ does indeed exist.) By \cite[Pr. 3.18]{K1}, the map $\fromto{\widetilde{\mathscr{F}}(\mathscr{C})}{N\Delta^{\op}}$ is a cocartesian fibration.

Dually, if $\mathscr{C}$ is a coWaldhausen $\infty$-category, write $\widetilde{\mathscr{F}}^{\vee}(\mathscr{C})$ for the simplicial set over $N\Delta$ satisfying the following universal property. We require, for any simplicial set $K$ and any map $\sigma\colon\fromto{K}{N\Delta}$, a bijection
\begin{equation*}
\Mor_{/N\Delta}(K,\widetilde{\mathscr{F}}^{\vee}(\mathscr{C}))\cong\Mor_{s\Set(2)}((K\times_{N\Delta}N\widetilde{\mathrm{M}}^{\op},K\times_{N\Delta}(N\widetilde{\mathrm{M}}^{\dag})^{\op}),(\mathscr{C},\mathscr{C}^{\dag})),
\end{equation*}
functorial in $\sigma$. By the dual of \cite[Pr. 3.18]{K1}, the map $\fromto{\widetilde{\mathscr{F}}^{\vee}(\mathscr{C})}{N\Delta}$ is a cartesian fibration, and it is clear that
\begin{equation*}
\widetilde{\mathscr{F}}^{\vee}(\mathscr{C})\cong\widetilde{\mathscr{F}}(\mathscr{C}^{\op})^{\op}.
\end{equation*}

The objects of the $\infty$-category $\widetilde{\mathscr{F}}(\mathscr{C})$ may be described as pairs $(m,X)$ consisting of a nonnegative integer $m$ and a functor $X\colon\fromto{\mathscr{O}(\Delta^m)}{\mathscr{C}}$ that carries ingressive morphisms to ingressive morphisms. Dually, the objects of the $\infty$-category $\widetilde{\mathscr{F}}^{\vee}(\mathscr{C})$ may be described as pairs $(m,X)$ consisting of a nonnegative integer $m$ and a functor $X\colon\fromto{\mathscr{O}(\Delta^m)^{\op}}{\mathscr{C}}$ that carries egressive morphisms to egressive morphisms.

Now if $\mathscr{C}$ is a Waldhausen $\infty$-category, we let $\widetilde{\mathscr{S}}(\mathscr{C})\subset\widetilde{\mathscr{F}}(\mathscr{C})$ denote the full subcategory spanned by those pairs $(m,X)$ such that for any integer $0\leq i\leq m$, the object $X(i,i)$ is a zero object of $\mathscr{C}$, and for any integers $0\leq i\leq k\leq j\leq\ell\leq m$, the square
\begin{equation*}
\begin{tikzpicture} 
\matrix(m)[matrix of math nodes, 
row sep=4ex, column sep=4ex, 
text height=1.5ex, text depth=0.25ex] 
{X(i,j)&X(i,\ell)\\ 
X(k,j)&X(k,\ell)\\}; 
\path[>=stealth,->,font=\scriptsize] 
(m-1-1) edge[>->] (m-1-2) 
edge[->] (m-2-1) 
(m-1-2) edge[->] (m-2-2) 
(m-2-1) edge[>->] (m-2-2); 
\end{tikzpicture}
\end{equation*}
is a pushout. Dually, if $\mathscr{C}$ is a coWaldhausen $\infty$-category, let $\widetilde{\mathscr{S}}^{\vee}(\mathscr{C})\subset\widetilde{\mathscr{F}}^{\vee}(\mathscr{C})$ denote the full subcategory spanned by those pairs $(m,X)$ such that for any integer $0\leq i\leq m$, the object $X(i,i)$ is a zero object of $\mathscr{C}$, and for any integers $0\leq i\leq k\leq j\leq\ell\leq m$, the square
\begin{equation*}
\begin{tikzpicture} 
\matrix(m)[matrix of math nodes, 
row sep=4ex, column sep=4ex, 
text height=1.5ex, text depth=0.25ex] 
{X(k,\ell)&X(k,j)\\ 
X(i,\ell)&X(i,j)\\}; 
\path[>=stealth,->,font=\scriptsize] 
(m-1-1) edge[->] (m-1-2) 
edge[->>] (m-2-1) 
(m-1-2) edge[->>] (m-2-2) 
(m-2-1) edge[->] (m-2-2); 
\end{tikzpicture}
\end{equation*}
is a pullback. Since ambigressive pullbacks and ambigressive pushouts coincide, we deduce that
\begin{equation*}
\widetilde{\mathscr{S}}^{\vee}(\mathscr{C})\cong\widetilde{\mathscr{S}}(\mathscr{C}^{\op})^{\op}.
\end{equation*}
\end{cnstr}

\begin{ntn} As in \cite{K1}, the constructions above yield functors
\begin{equation*}
\widetilde{\mathscr{S}}\colon\fromto{\Wald_{\infty}}{\Cat_{\infty,/N\Delta^{\op}}}\textrm{\quad and\quad}\widetilde{\mathscr{S}}^{\vee}\colon\fromto{\coWald_{\infty}}{\Cat_{\infty,/N\Delta}}.
\end{equation*}
The functor $\fromto{\mathrm{M}}{\widetilde{\mathrm{M}}}$ given by the assignment $\goesto{(m,i)}{(m,0,i)}$ induces a natural transformation $\fromto{\widetilde{\mathscr{S}}}{\mathscr{S}}$ over $N\Delta^{\op}$ and a natural transformation $\fromto{\widetilde{\mathscr{S}}^{\vee}}{\mathscr{S}^{\vee}}$ over $N\Delta$.

Furthermore, we can pass to the functors that classify these fibrations to obtain functors
\begin{equation*}
\widetilde{\SS}_{\ast}\colon\fromto{\Wald_{\infty}}{\Fun(N\Delta^{\op},\Cat_{\infty})}\textrm{\quad and\quad}\widetilde{\SS}_{\ast}^{\vee}\colon\fromto{\coWald_{\infty}}{\Fun(N\Delta^{\op},\Cat_{\infty})}.
\end{equation*}
For any Waldhausen $\infty$-category (respectively, any coWaldhausen $\infty$-category) $\mathscr{C}$, the simplicial category $\widetilde{\SS}_{\ast}(\mathscr{C})\colon\fromto{N\Delta^{\op}}{\Cat_{\infty}}$ (resp., $\widetilde{\SS}^{\vee}_{\ast}(\mathscr{C})\colon\fromto{N\Delta^{\op}}{\Cat_{\infty}}$) carries an object $\mathbf{m}$ to the full subcategory
\begin{equation*}
\widetilde{\SS}_m(\mathscr{C})\subset\Fun(\mathscr{O}(\Delta^m),\mathscr{C})\textrm{\qquad (resp.,\quad}\widetilde{\SS}^{\vee}_m(\mathscr{C})\subset\Fun(\mathscr{O}(\Delta^m)^{\op},\mathscr{C})\textrm{\quad).}
\end{equation*}
spanned by those diagrams $X$ such that for any integer $0\leq i\leq m$, the object $X(i,i)$ is a zero object of $\mathscr{C}$, and for any integers $0\leq i\leq k\leq j\leq\ell\leq m$, the square
\begin{equation*}
\begin{tikzpicture} 
\matrix(m)[matrix of math nodes, 
row sep=4ex, column sep=4ex, 
text height=1.5ex, text depth=0.25ex] 
{X(i,j)&X(i,\ell)\\ 
X(k,j)&X(k,\ell)\\}; 
\path[>=stealth,->,font=\scriptsize] 
(m-1-1) edge[>->] (m-1-2) 
edge (m-2-1) 
(m-1-2) edge (m-2-2) 
(m-2-1) edge[>->] (m-2-2); 
\end{tikzpicture}
\end{equation*}
is an ambigressive pushout (resp., for any integers $0\leq i\leq k\leq j\leq\ell\leq m$, the square
\begin{equation*}
\begin{tikzpicture} 
\matrix(m)[matrix of math nodes, 
row sep=4ex, column sep=4ex, 
text height=1.5ex, text depth=0.25ex] 
{X(k,\ell)&X(k,j)\\ 
X(i,\ell)&X(i,j)\\}; 
\path[>=stealth,->,font=\scriptsize] 
(m-1-1) edge (m-1-2) 
edge[->>] (m-2-1) 
(m-1-2) edge[->>] (m-2-2) 
(m-2-1) edge (m-2-2); 
\end{tikzpicture}
\end{equation*}
is an ambigressive pullback).
\end{ntn}

In light of the uniqueness of limits and colimits in $\infty$-categories \cite[Pr. 1.2.12.9]{HTT}, one readily has the following.
\begin{prp} Suppose $\mathscr{C}$ a Waldhausen $\infty$-category. Then the functor
\begin{equation*}
\fromto{\widetilde{\mathscr{S}}(\mathscr{C})}{N\Delta^{\op}}
\end{equation*}
is a cocartesian fibration, and the map $\fromto{\widetilde{\mathscr{S}}(\mathscr{C})}{\mathscr{S}(\mathscr{C})}$ defined above is a fiberwise equivalence over $N\Delta^{\op}$. Dually, if $\mathscr{C}$ is a coWaldhausen $\infty$-category, then the functor
\begin{equation*}
\fromto{\widetilde{\mathscr{S}}^{\vee}(\mathscr{C})}{N\Delta}
\end{equation*}
is a cartesian fibration, and the map $\fromto{\widetilde{\mathscr{S}}^{\vee}(\mathscr{C})}{\mathscr{S}^{\vee}(\mathscr{C})}$ defined above is a fiberwise equivalence over $N\Delta$.
\end{prp}

\begin{nul} In particular, we can use the equivalence of the previous proposition to endow, for any Waldhausen $\infty$-category (respectively, any coWaldhausen $\infty$-category) $\mathscr{C}$, the $\infty$-category $\widetilde{\mathscr{S}}(\mathscr{C})$ (resp., $\widetilde{\mathscr{S}}^{\vee}(\mathscr{C})$) with the structure of a Waldhausen $\infty$-category. That is, let us declare a morphism $\fromto{X}{Y}$ of $\widetilde{\mathscr{S}}(\mathscr{C})$ (resp., of $\widetilde{\mathscr{S}}^{\vee}(\mathscr{C})$) ingressive (resp., egressive) just in case its image in $\mathscr{S}(\mathscr{C})$ (resp., in in $\mathscr{S}^{\vee}(\mathscr{C})$) is so. Observe that under this definition, $\fromto{X}{Y}$ is ingressive (resp., egressive) just in case each of the morphisms $\fromto{X(i,j)}{Y(i,j)}$ of $\mathscr{C}$ is so.

Consequently, the constructions $\widetilde{\SS}$ and $\widetilde{\SS}^{\vee}$ yield functors
\begin{equation*}
\widetilde{\SS}\colon\fromto{\Wald_{\infty}}{\Fun(N\Delta^{\op},\Wald_{\infty})}\textrm{\quad and\quad}\widetilde{\SS}^{\vee}\colon\fromto{\coWald_{\infty}}{\Fun(N\Delta^{\op},\coWald)}.
\end{equation*}
\end{nul}

Now suppose $\mathscr{C}$ an exact $\infty$-category. Since ambigressive pushouts and ambigressive pullbacks coincide in $\mathscr{C}$, it follows that there is a canonical equivalence
\begin{equation*}
(\widetilde{\SS}_{\ast}\circ\op)|_{\Exact_{\infty}}\simeq(\widetilde{\SS}^{\vee}_{\ast})|_{\Exact_{\infty}},
\end{equation*}
where $\op\colon\equivto{N\Delta^{\op}}{N\Delta^{\op}}$ is the opposite automorphism of $N\Delta^{\op}$. We therefore deduce the following.
\begin{thm} Suppose $\mathscr{E}$ an $\infty$-topos. For any pointed functor $\rho\colon\fromto{\Cat_{\infty}^{\ast}}{\mathscr{E}_{\ast}}$ that preserves filtered colimits, an equivalence $\equivto{\rho}{\rho\circ\op}$ induces a canonical exact duality on the Goodwillie additivization $D\rho$.
\begin{proof} The equivalence $\equivto{\rho}{\rho\circ\op}$, combined with the equivalence $\widetilde{\SS}_{\ast}\circ\op\simeq\widetilde{\SS}^{\vee}_{\ast}$, yields an equivalence
\begin{equation*}
\equivto{|\rho\circ\widetilde{\SS}_{\ast}|\simeq|\rho\circ\widetilde{\SS}_{\ast}\circ\op|}{|\rho\circ\op\circ\widetilde{\SS}_{\ast}\circ\op|\simeq|\rho^{\vee}\circ\widetilde{\SS}^{\vee}_{\ast}|}.\qedhere
\end{equation*}
\end{proof}
\end{thm}
\noindent Applying this result to the functor $\iota$ yields the following.
\begin{cor}\label{cor:Kdual} Algebraic $K$-theory admits an exact duality,
\end{cor}

\section{Theorem of the heart} In this section, we show that the Waldhausen $K$-theory of a stable $\infty$-category with a bounded t-structure agrees with the $K$-theory of its heart. Amnon Neeman has provided an analogous result for his $K$-theory of triangulated categories \cite{MR1930883}; given Neeman's result, our result here may be alternatively summarized as saying that the Waldhausen $K$-theory of a stable $\infty$-category $\mathscr{A}$ with a bounded t-structure agrees with Neeman's $K$-theory of the triangulated homotopy category $\mathscr{T}=h\mathscr{A}$ (denoted $K({}^w\!\mathscr{T})$ in \cite{MR2181838}); this verifies a conjecture posed by Neeman \cite[App. A]{MR1793672} for such stable homotopy theories.

In this section, suppose $E$ a small stable $\infty$-category equipped with a bounded t-structure $(E_{\leq-1},E_{\geq 0})$. Our objective is to prove the following
\begin{thm}[Heart]\label{thm:thmofheart} The inclusions $\into{E^{\heartsuit}}{E_{\geq 0}}$ and $\into{E_{\geq 0}}{E}$ induce equivalences
\begin{equation*}
\equivto{K(E^{\heartsuit})}{K(E_{\geq 0})}\textrm{\quad and\quad}\equivto{K(E_{\geq 0})}{K(E)}.
\end{equation*}
\end{thm}
The fact that the inclusions induce isomorphisms $K_0(E^{\heartsuit})\cong K_0(E_{\geq0})\cong K_0(E)$ is well known and trivial. Consequently, appealing to the Cofinality Theorem of \cite[Th. 10.11]{K1} (which states that idempotent completions induce equivalences on the connected cover of $K$-theory), we may therefore assume that $E$ is idempotent complete. We now set about proving that the higher $K$-groups of $E^{\heartsuit}$ and $E$ agree.

Our proof is quite straightforward. The first main tool is the following.
\begin{thm}[Special Fibration Theorem, \protect{\cite[Pr 10.12]{K1}}]\label{cor:specialfibration} Suppose $C$ a compactly generated $\infty$-category containing a zero object, suppose $L\colon\fromto{C}{D}$ an accessible localization, and suppose the inclusion $\into{D}{C}$ preserves filtered colimits. Assume also that the class of all $L$-equivalences of $C$ is generated (as a strongly saturated class) by the $L$-equivalences between compact objects. Then $L$ induces a pullback square of spaces
\begin{equation*}
\begin{tikzpicture} 
\matrix(m)[matrix of math nodes, 
row sep=4ex, column sep=4ex, 
text height=1.5ex, text depth=0.25ex] 
{K(E^{\omega})&K(C^{\omega})\\ 
\ast&K(D^{\omega}),\\}; 
\path[>=stealth,->,font=\scriptsize] 
(m-1-1) edge (m-1-2) 
edge (m-2-1) 
(m-1-2) edge (m-2-2) 
(m-2-1) edge (m-2-2); 
\end{tikzpicture}
\end{equation*}
where $C^{\omega}$ and $D^{\omega}$ are equipped with the maximal pair structure, and $E^{\omega}\subset C^{\omega}$ is the full subcategory spanned by those objects $X$ such that $LX\simeq 0$.
\end{thm}

As a consequence of \cite[Cor. 5.5.7.4 and Pr. 5.5.7.8]{HTT}, we proved in \cite[Cor 10.12.2]{K1} that the functor $\tau_{\leq -1}$ induces a pullback square
\begin{equation*}
\begin{tikzpicture} 
\matrix(m)[matrix of math nodes, 
row sep=4ex, column sep=4ex, 
text height=1.5ex, text depth=0.25ex] 
{K(E_{\geq 0})&K(E)\\ 
\ast&K(E_{\leq -1}^{\mathit{max}}),\\}; 
\path[>=stealth,->,font=\scriptsize] 
(m-1-1) edge (m-1-2) 
edge (m-2-1) 
(m-1-2) edge (m-2-2) 
(m-2-1) edge (m-2-2); 
\end{tikzpicture}
\end{equation*}
of spaces, where the $\infty$-categories that appear are equipped with their maximal Waldhausen structure --- i.e., the one in which every morphism is ingressive. (For $E_{\geq0}$ and $E$, the Waldhausen structure described in Ex. \ref{exm:tstructexacts} \emph{is} the maximal Waldhausen structure, but this is very much not the case for $E_{\leq -1}$.) 

We now claim that the $K$-theory of the maximal Waldhausen $\infty$-category $E_{\leq -1}^{\mathit{max}}$. Indeed, we may apply \cite[Cor. 8.2.1]{K1} to the $\infty$-category $E_{\leq -1}$; this will ensure that the functor
\begin{equation*}
\Sigma^{\infty}_{E_{\leq -1}}\colon\fromto{E_{\leq -1}^{\mathit{max}}}{\widetilde{\Sp}\left(E_{\leq -1}\right)^{\mathit{max}}}
\end{equation*}
induces an equivalence on $K$-theory. Note that the suspension functor on $E_{\leq -1}$ is the composite $\tau_{\leq -1}\circ\Sigma_E$. Since the t-structure is bounded, it therefore follows that $\Sigma^{\infty}_{E_{\leq -1}}$ is equivalent to the constant functor at $0$.

\begin{wrn} Note that this argument applies only to the maximal Waldhausen $\infty$-category $E_{\leq -1}^{\mathit{max}}$. The $K$-theory of the Waldhausen structure on $E_{\leq -1}$ of Ex. \ref{exm:tstructexacts} --- in which a morphism is ingressive just in case it induces a monomorphism on $\pi_{-1}$ --- will turn out to agree with the $K$-theory of $E$.
\end{wrn}

Consequently, the inclusion $\into{E_{\geq0}}{E}$ induces an equivalence
\begin{equation*}
\equivto{K\left(E_{\geq0}\right)}{K(E)},
\end{equation*}
and we are left with showing that the map $\fromto{K(E^{\heartsuit})}{K(E_{\geq0})}$ is an equivalence.

For this, we take opposites. The opposite $\infty$-category $E^{\op}$ is endowed with the dual t-structure, with $\left(E^{\op}\right)_{\leq-n}=\left(E_{\geq n}\right)^{\op}$. Consider the functor
\begin{equation*}
\tau_{\leq-1}=(\tau_{\geq1})^{\op}\colon\fromto{(E_{\geq 0})^{\op}=(E^{\op})_{\leq0}}{(E^{\op})_{\leq-1}=(E_{\geq1})^{\op}}.
\end{equation*}
Our claim is that it induces a pullback square
\begin{equation}\label{finalnail}\tag{$\ast$}
\begin{tikzpicture}[baseline]
\matrix(m)[matrix of math nodes, 
row sep=4ex, column sep=4ex, 
text height=1.5ex, text depth=0.25ex] 
{K((E^{\op})^{\heartsuit})&K((E^{\op})_{\leq0})\\ 
0&K((E^{\op})^{\mathit{max}}_{\leq-1}),\\}; 
\path[>=stealth,->,font=\scriptsize] 
(m-1-1) edge (m-1-2) 
edge (m-2-1) 
(m-1-2) edge (m-2-2) 
(m-2-1) edge (m-2-2); 
\end{tikzpicture}
\end{equation}
where the $\infty$-categories that appear are equipped with the Waldhausen structure in which a morphism is ingressive just in case it induces a monomorphism (in $E^{\op}$) on $\pi_0$ (and in particular, the Waldhausen structure on $(E^{\op})_{\leq-1}$ is maximal).

Of course $K((E^{\op})^{\mathit{max}}_{\leq-1})$ vanishes just as $K(E^{\mathit{max}}_{\leq -1})$ vanishes. Consequently, if we verify that \eqref{finalnail} is a pullback, we will deduce that the map $\fromto{K((E^{\op})^{\heartsuit})}{K((E^{\op})_{\leq0})}$ is an equivalence, and in light of Cor. \ref{cor:Kdual}, which allows us to pass between a Waldhausen $\infty$-category and its opposite under $K$-theory, the proof of Th. \ref{thm:thmofheart} will be complete.

To prove that \eqref{finalnail} is a pullback, we cannot simply appeal to the Special Fibration Theorem, because the Waldhausen $\infty$-categories that appear have non-maximal pair structures. Instead, we must appeal to what we have called the Generic Fibration Theorem II; this is the $\infty$-categorical variant of Waldhausen's celebrated Fibration Theorem. Unfortunately, this means that we will have to check some technical hypotheses, but there is no escape.
\begin{thm}[Generic Fibration Theorem II, \protect{\cite[Th. 9.24]{K1}}]\label{genfib} Suppose $C$ a Waldhausen $\infty$-category, and suppose that $wC$ is a subcategory of $C$ that satisfies the following conditions.
\begin{enumerate}[(\ref{genfib}.1)]
\item Every equivalence of $C$ lies in $wC$.
\item The morphisms of $wC$ satisfy (the $\infty$-categorical analogue of) \emph{Waldhausen's gluing axiom}. That is, for any cofibrations $\cofto{U}{V}$ and $\cofto{X}{Y}$ and any cube
\begin{equation*}
\begin{tikzpicture}[cross line/.style={preaction={draw=white, -, 
line width=6pt}}]
\matrix(m)[matrix of math nodes, 
row sep=2ex, column sep=2ex, 
text height=1.5ex, text depth=0.25ex]
{&U&&V\\
U'&&V'&\\
&X&&Y\\
X'&&Y'&\\
};
\path[>=stealth,->,font=\scriptsize]
(m-1-2) edge (m-2-1)
edge (m-3-2)
edge[>->] (m-1-4)
(m-3-2) edge (m-4-1)
edge[>->] (m-3-4)
(m-2-1) edge[cross line,>->] (m-2-3)
edge (m-4-1)
(m-1-4) edge (m-2-3)
edge (m-3-4)
(m-4-1) edge[>->] (m-4-3)
(m-3-4) edge (m-4-3)
(m-2-3) edge[cross line] (m-4-3);
\end{tikzpicture}
\end{equation*}
in which the top and bottom faces are pushout squares, if $\fromto{U}{X}$, $\fromto{V}{Y}$ and $\fromto{U'}{X'}$ all lie in $wC$, then so does $\fromto{V'}{Y'}$.
\item There are \emph{enough cofibrations} in the following sense. The three functors
\begin{equation*}
\into{wC\cap C_{\dag}}{wC},\textrm{\qquad}\into{w\Fun(\Delta^1,C)\cap\Fun(\Delta^1,C)_{\dag}}{w\Fun(\Delta^1,C)},
\end{equation*}
and
\begin{equation*}
\into{w\FF_1(C)\cap\FF_1(C)_{\dag}}{w\FF_1(C)}
\end{equation*}
are all weak homotopy equivalences, where cofibrations of $\Fun(\Delta^1,C)$ are defined objectwise, and $w\Fun(\Delta^1,C)$ and $w\FF_1(C)$ are also defined objectwise.
\end{enumerate}
Now denote by $C^w\subset C$ the full subcategory spanned by those objects $X$ such that the cofibration $\fromto{0}{X}$ lies in $wC$, and declare a morphism therein to be ingressive just in case it is so in $C$. Denote by $\BB_{\ast}(C,wC)$ the simplicial $\infty$-category whose $\infty$-category of $m$-simplices is the full subcategory $\BB_m(C,wC)\subset\Fun(\Delta^m,C)$ spanned by those sequences of edges 
\begin{equation*}
X_0\ \tikz[baseline]\draw[>=stealth,->,font=\scriptsize,inner sep=0.5pt](0,0.5ex)--(0.5,0.5ex);\ X_1\ \tikz[baseline]\draw[>=stealth,->,font=\scriptsize,inner sep=0.5pt](0,0.5ex)--(0.5,0.5ex);\ \cdots\ \tikz[baseline]\draw[>=stealth,->,font=\scriptsize,inner sep=0.5pt](0,0.5ex)--(0.5,0.5ex);\ X_m
\end{equation*}
with the property that each $\fromto{X_i}{X_j}$ lies in $wC$; declare a morphism of $\BB_m(C,wC)$ to be a cofibration just in case it is so objectwise. Then $C^w$ is a Waldhausen $\infty$-category, and $\BB_{\ast}(C,wC)$ is a simplicial Waldhausen $\infty$-category, and the obvious functors induce a fiber sequence of spaces
\begin{equation*}
\begin{tikzpicture} 
\matrix(m)[matrix of math nodes, 
row sep=4ex, column sep=4ex, 
text height=1.5ex, text depth=0.25ex] 
{K(C^w)&K(C)\\ 
\ast&\left|K(\BB_{\ast}(C,wC))\right|.\\}; 
\path[>=stealth,->,font=\scriptsize] 
(m-1-1) edge (m-1-2) 
edge (m-2-1) 
(m-1-2) edge (m-2-2) 
(m-2-1) edge (m-2-2); 
\end{tikzpicture}
\end{equation*}
Here of course $|\cdot|$ denotes the geometric realization.
\end{thm}

We apply this theorem to the Waldhausen $\infty$-category $(E^{\op})_{\leq0}$ along with the subcategory $w((E^{\op})_{\leq0})$ consisting of those morphisms $\fromto{X}{Y}$ such that the induced morphism $\fromto{\tau_{\leq-1}X}{\tau_{\leq-1}Y}$ is an equivalence. It is an easy matter to check the gluing axiom in this setting. To show that we have enough cofibrations, let us employ the following construction: factor any morphism $f\colon\fromto{X}{Y}$ of $wC$ as
\begin{equation*}
X\ \tikz[baseline]\draw[>=stealth,->,font=\scriptsize,inner sep=0.5pt](0,0.5ex)--(0.5,0.5ex);\ X_f\ \tikz[baseline]\draw[>=stealth,>->,font=\scriptsize,inner sep=0.5pt](0,0.5ex)--(0.5,0.5ex);\ Y,
\end{equation*}
where $X_f$ is the fiber of the natural map $\fromto{Y}{\tau_{\leq 0}(Y/X)}$, where $Y/X$ is the cofiber of $f$. This construction defines a deformation retraction $\fromto{\Fun(\Delta^1,wC)}{\Fun(\Delta^1,wC\cap C_{\dag})}$ of the inclusion $\into{\Fun(\Delta^1,wC\cap C_{\dag})}{\Fun(\Delta^1,wC)}$. It follows from the functoriality of this construction that it also defines a deformation retraction of the inclusions
\begin{equation*}
\into{\Fun(\Delta^1,w\Fun(\Delta^1,C)\cap\Fun(\Delta^1,C)_{\dag})}{\Fun(\Delta^1,w\Fun(\Delta^1,C))},
\end{equation*}
and
\begin{equation*}
\into{\Fun(\Delta^1,w\FF_1(C)\cap\FF_1(C)_{\dag})}{\Fun(\Delta^1,w\FF_1(C))}.
\end{equation*}
Consequently, these maps are all homotopy equivalences. 

So Th. \ref{genfib} now ensures that we have a pullback square
\begin{equation*}
\begin{tikzpicture} 
\matrix(m)[matrix of math nodes, 
row sep=4ex, column sep=4ex, 
text height=1.5ex, text depth=0.25ex] 
{K((E^{\op})^{\heartsuit})&K((E^{\op})_{\leq0})\\ 
0&\left|K(\BB_{\ast}((E^{\op})_{\leq0},w(E^{\op})_{\leq0}))\right|,\\}; 
\path[>=stealth,->,font=\scriptsize] 
(m-1-1) edge (m-1-2) 
edge (m-2-1) 
(m-1-2) edge (m-2-2) 
(m-2-1) edge (m-2-2); 
\end{tikzpicture}
\end{equation*}
and it remains only to identify the cofiber term. Now it follows from \cite[Pr. 10.10]{K1} that, in the situation of Th. \ref{genfib}, the geometric realization $\left|K(\BB_{\ast}(C,wC))\right|$ can be identified with $\Omega|w\SS_{\ast}(C)|$, where $w\SS_m(C)\subset\SS_m(C)$ is the subcategory whose morphisms are weak equivalences. Consequently, to identify $K((E^{\op})^{\mathit{max}}_{\leq-1})$ with the geometric realization
\begin{equation*}
\left|K(\BB_{\ast}((E^{\op})_{\leq0},w(E^{\op})_{\leq0}))\right|,
\end{equation*}
it is enough to show that the functor
\begin{equation*}
p\colon\fromto{w\SS_m((E^{\op})_{\leq0})}{\iota\SS_m((E^{\op})^{\mathit{max}}_{\leq-1})}
\end{equation*}
induced by $\tau_{\leq-1}$ is a weak homotopy equivalence. Let us show that the inclusion
\begin{equation*}
i\colon\into{\iota\SS_m((E^{\op})^{\mathit{max}}_{\leq-1})}{\iota\SS_m((E^{\op})_{\leq0})\subset w\SS_m((E^{\op})_{\leq0})}
\end{equation*}
defines a homotopy inverse. Since $\tau_{\leq-1}$ is a localization functor, we obtain for every $m\geq 0$ a natural equivalence $p\circ i\simeq\id$. In the other direction, the natural transformation $\fromto{\id}{\tau_{\leq-1}}$ induces a natural transformation $\fromto{\id}{i\circ p}$. Hence $i$ and $p$ are homotopy inverse.\qed

Amnon Neeman's Theorem of the Heart now implies the following, which verifies some cases of his conjecture \cite[Conj. A.5]{MR1793672}.
\begin{cor}\label{cor:Neemanconj} For any idempotent-complete stable $\infty$-category $\mathscr{A}$, if the triangulated homotopy category $\mathscr{T}=h\mathscr{A}$ admits a bounded t-structure, then we have canonical equivalences
\begin{equation*}
K(\mathscr{A})\simeq K(\mathscr{T}^{\heartsuit})\simeq K({}^w\!\mathscr{T}).
\end{equation*}
\end{cor}


\section{Application: Abelian models for the algebraic $G$-theory of schemes} A trivial application of the Theorem of the Heart applies to is that the $K$-theory of an abelian category $A$ agrees with the $K$-theory of its bounded derived $\infty$-category $D^b(A)$ \cite[1.11.7]{MR92f:19001}. However, tilting theory provides other bounded t-structures on the $\infty$-category $\mathrm{D}^{b}(A)$. The $K$-theory of the heart of any of these t-structures will agree with the $K$-theory of $A$. Let us explore one class of examples now.

\begin{exm}\label{exm:perv} Suppose $X$ a noetherian scheme equipped with a dualizing complex $\omega_X\in\mathrm{D}^{b}(\Coh(X))$. Then Arinkin and Bezrukavnikov \cite{MR2668828}, following Deligne, construct a family of t-structures on $\mathrm{D}^{b}(\Coh(X))$ in the following manner. (Here we use \emph{cohomological} grading conventions, to maintain compatibility with \cite{MR2668828}.) Write $X^{\mathrm{top}}$ for the for the underlying topological space of $X$, and define $\dim\colon\fromto{X^{\mathrm{top}}}{\ZZ}$ a map such that $i_x^{!}\omega_X$ is concentrated in degree $-\dim(x)$. Suppose $p\colon\fromto{X^{\mathrm{top}}}{\ZZ}$ a function --- called a \emph{perversity} --- such that for any points $z,x\in X^{\mathrm{top}}$ such that $z\in\overline{\{x\}}$, one has
\begin{equation*}
p(x)\leq p(z)\leq p(x)+\dim(x)-\dim(z).
\end{equation*}
Let $\mathrm{D}^{p\geq 0}\subset\mathrm{D}^b(\Coh(X))$ be the full subcategory spanned by those complexes $F$ such that for any point $x\in X^{\mathrm{top}}$, one has $i_x^!F\in\mathrm{D}^{\geq p(x)}(\mathscr{O}_{X,x})$; dually, let $\mathrm{D}^{p\leq 0}\subset\mathrm{D}^b(\Coh(X))$ be the full subcategory spanned by those complexes $F$ such that for any point $x\in X^{\mathrm{top}}$, one has $i_{x,\star}F\in\mathrm{D}^{\leq p(x)}(\mathscr{O}_{X,x})$. Then $(\mathrm{D}^{p\geq 0},\mathrm{D}^{p\leq 0})$ define a bounded t-structure on $\mathrm{D}^{b}(\Coh(X))$ \cite[Th. 3.10]{MR2668828}.

The algberaic $K$-theory of the heart $\mathrm{D}^{p,\heartsuit}$ of this t-structure now coincides with the $G$-theory of $X$. Let us list two special cases of this.
\begin{enumerate}[(\ref{exm:perv}.1)]
\item Suppose $S$ a set of prime numbers. Let $\mathscr{E}_S$ be the full subcategory of $\mathrm{D}^b(\Coh(\ZZ))$ generated under extensions by the objects
\begin{equation*}
\ZZ,\qquad\{\ZZ/p\ |\ p\in S\},\qquad\{\ZZ/p[1]\quad |\quad p\notin S\}.
\end{equation*}
Then $\mathscr{E}_S$ is an abelian category whose $K$-theory coincides with the $K$-theory of $\ZZ$.
\item For any noetherian scheme with a dualizing complex $\omega_X\in\mathrm{D}^{b}(\Coh(X))$, the $K$-theory of the abelian category of \emph{Cohen--Macaulay complexes} (i.e., those complexes $F\in\mathrm{D}^{b}(\Coh(X))$ such that the complex
\begin{equation*}
\mathbf{D}F\coloneq\RR\Mor_{\mathscr{O}_X}(F,\omega_X)
\end{equation*}
is concentrated in degree $0$, \cite[\S 6]{MR2244264}) agrees with the $G$-theory of $X$.
\end{enumerate}
\end{exm}


\section{Application: A theorem of Blumberg--Mandell} In this section, we give a new proof of the theorem of Blumberg--Mandell \cite{BM} that establishes a localization sequence
\begin{equation*}
K(\pi_0E)\to K(e)\to K(E)
\end{equation*}
for any suitable even periodic $E_1$ ring spectrum $E$ with $\pi_0E$ regular noetherian, where $e$ denotes the connective cover of $E$. In light of \cite[Pr. 13.16]{K1}, the key point is the identification of the fiber term; this is the subject of this section. 

Recall \cite[Pr. 8.2.5.16]{HA} that a connective $E_1$ ring $\Lambda$ is said to be \emph{left coherent} if $\pi_0\Lambda$ is left coherent as an ordinary ring, and if for any $n\geq 1$, the left $\pi_0\Lambda$-module $\pi_n\Lambda$ is finitely presented.

A left module $M$ over a left coherent $E_1$ ring $\Lambda$ is \emph{almost perfect} just in case $\pi_mM=0$ for $m\ll 0$ and for any $n$, the left $\pi_0\Lambda$-module $\pi_nM$ is finitely presented \cite[Pr. 8.2.5.17]{HA}.

\begin{dfn} Suppose $\Lambda$ a left coherent $E_1$ ring, and suppose $M$ a left $\Lambda$-module. We say that $M$ is \emph{truncated} if $\pi_mM=0$ for $m\gg 0$. We say that $M$ is \emph{coherent} if it is almost perfect and truncated. Write $\Coh^{\ell}_{\Lambda}\subset\Mod^{\ell}_{\Lambda}$ for the full subcategory spanned by the coherent modules, and write $G(\Lambda)$ for $K(\Coh^{\ell}_{\Lambda})$.
\end{dfn}

\begin{wrn} In general, it is not the case that a perfect $\Lambda$-module is coherent; consequently, the usual Poincar\'e duality map $\fromto{K}{G}$ for discrete rings does not have an obvious analogue for $E_1$ rings.
\end{wrn}

It turns out that $G$-theory is not a new invariant of $E_1$ rings, since we have the following new proof of the D\'evissage Theorem of Blumberg--Mandell \cite{BM}.
\begin{prp}\label{prp:dev} For any coherent $E_1$ ring $\Lambda$, the inclusion $\into{N\Mod^{\ell,\mathrm{fp}}_{\pi_0\Lambda}}{\Coh^{\ell}_{\Lambda}}$ of the nerve of the category of finitely presented $\pi_0\Lambda$-modules induces an equivalence
\begin{equation*}
\equivto{G(\pi_0\Lambda)}{G(\Lambda)}.
\end{equation*}
\begin{proof} We note that $\Coh^{\ell}_{\Lambda}$ is the full subcategory of the $\infty$-category of almost perfect $\Lambda$-modules spanned by those that are bounded for the t-structure given by \cite[Pr. 8.2.5.18]{HA}. Furthermore, \cite[Pr. 8.2.5.11(2)]{HA} applies to ensure that $\Coh^{\ell}_{\Lambda}$ is idempotent complete. Consequently, the Theorem of the Heart applies, and the proof is complete once one observes that the heart $\Coh^{\ell,\heartsuit}_{\Lambda}$ can be identified with $N\Mod^{\ell,\mathrm{fp}}_{\pi_0\Lambda}$ \cite[Rk. 8.2.5.19]{HA}.
\end{proof}
\end{prp}
\noindent Consequently, from the point of view of ``brave new algebra,'' $G$-theory is a relatively coarse invariant.

Now we hope to compare the $G$-theory of an $E_1$ ring to the $K$-theory of the $\infty$-category of truncated perfect modules. This requires a weak regularity hypothesis, which we formulate now.

\begin{dfn} Let us say that a coherent $E_1$ ring $\Lambda$ is \emph{almost regular} if any truncated, almost perfect $\Lambda$-module has finite Tor dimension.
\end{dfn}

\begin{exm} If the graded ring $\pi_{\ast}\Lambda$ has finite Tor-dimension (e.g., if $\pi_{\ast}\Lambda$ is a regular noetherian ring), then the Tor spectral sequence ensures that $\Lambda$ is almost regular.
\end{exm}

The following result is now an immediate consequence of \cite[Pr. 8.2.5.23(4)]{HA}
\begin{prp} Suppose $\Lambda$ a coherent $E_1$ ring that is almost regular. Then the $\infty$-category $\Perf^{\ell,b}_{\Lambda}$ of perfect truncated left $\Lambda$-modules coincides with the $\infty$-category $\Coh^{\ell}_{\Lambda}$ of coherent left $\Lambda$-modules.
\end{prp}

Assembling all this, we obtain the following formulation of the Localization Theorem of \cite{BM}.
\begin{thm}\label{thm:BM} Suppose $\Lambda$ a coherent $E_1$ ring spectrum that is almost regular, and suppose $L\colon\fromto{\Mod^{\ell}(\Lambda)}{\Mod^{\ell}(\Lambda)}$ smashing Bousfield localization with the property that a left $\Lambda$-module $M$ is $L$-acyclic just in case it is truncated. (Note that in this case, $L$ is automatically a finite localization functor.) Then there is a fiber sequence of spaces
\begin{equation*}
G(\pi_0\Lambda)\to K(\Lambda)\to K(L(\Lambda)),
\end{equation*}
\end{thm}

Of course when $\pi_0\Lambda$ is regular, the fiber term can be identified with $K(\pi_0\Lambda)$.

\begin{exm}\label{exm:BMfiberseqs} Here are some examples of fiber sequences resulting from this theorem.
\begin{enumerate}[(\ref{exm:BMfiberseqs}.1)]
\item Consider the Adams summand $L$ with its canonical $E_{\infty}$ structure; its connective cover $\ell$ admits a canonical $E_{\infty}$ as well \cite{MR2358512}. The fiber sequence above becomes
\begin{equation*}
K(\ZZ)\to K(\ell)\to K(L).
\end{equation*}
\item Similarly, one can use the $E_{\infty}$ structure on complex $K$-theory $K\!U$ and on its connective cover to obtain
\begin{equation*}
K(\ZZ)\to K(ku)\to K(KU).
\end{equation*}
\item For any perfect field $k$ of characteristic $p>0$ and any formal group $\Gamma$ of height $n$ over $k$, consider the Lubin--Tate spectrum $E(k,\Gamma)$ with its canonical $E_{\infty}$ structure and its connective cover $e(k,\Gamma)$ with its induced $E_{\infty}$ structure \cite{MR2358512}. In this case, the fiber sequence above becomes
\begin{equation*}
K(\WW(k)[\![u_1,\dots,u_{n-1}]\!])\to K(e(k,\Gamma))\to K(E(k,\Gamma)).
\end{equation*}
\item Given any $E_1$ structure on Morava $K$-theory $K(n)$ and a compatible one on its connective cover $k(n)$, one has a fiber sequence
\begin{equation*}
K(\FF_p)\to K(k(n))\to K(K(n)).
\end{equation*}
\end{enumerate}
\end{exm}


\section{Application: $G$-theory of spectral Deligne--Mumford stacks} The purpose of this final very brief section is simply to note that the $G$-theory of locally noetherian spectral Deligne--Mumford stacks is insensitive to derived structure.

\begin{dfn} For any spectral Deligne--Mumford stack $\mathscr{X}$, we write $\Coh(\mathscr{X})$ for the stable $\infty$-category of coherent sheaves on $\mathscr{X}$ \cite[Df. 2.6.20]{DAGVIII}, i.e., those quasicoherent sheaves that are almost perfect and locally truncated. Write $G(\mathscr{X})$ for the algebraic $K$-theory of $\Coh(\mathscr{X})$.
\end{dfn}

Now the Theorem of the Heart, combined with \cite[Rk. 2.3.20]{DAGVIII}, instantly yields the following.

\begin{prp}\label{prp:Gthynothingnew} For any locally noetherian spectral Deligne--Mumford stack $\mathscr{X}$ with underlying ordinary Deligne--Mumford stack $\mathscr{X}_0$, the embedding
\begin{equation*}
\into{N\Coh(\mathscr{X}_0)}{\Coh(\mathscr{X})}
\end{equation*}
induces an equivalence
\begin{equation*}
\equivto{G(\mathscr{X}_0)}{G(\mathscr{X})}.
\end{equation*}
\end{prp}

\noindent Roughly speaking, just as $G$-theory is invariant under ordinary nilpotent thickenings, it turns out that $G$-theory is invariant under \emph{derived} nilpotent thickenings as well.


\bibliographystyle{amsplain}
\bibliography{kthy}

\end{document}